\documentclass{elsarticle}

\usepackage{amssymb,amsthm,amsmath}
\usepackage[usenames,dvipsnames]{xcolor}
\usepackage{soul}
\usepackage[normalem]{ulem}
\usepackage{cancel}
\usepackage{mathabx}
\allowdisplaybreaks
\usepackage{color,graphicx,epsfig}
\usepackage{mathrsfs}
\usepackage{enumerate} 
\usepackage{enumitem} 
\usepackage{esvect}

\newcommand{\comment}[1]{}

\definecolor{teal}{RGB}{0,128,128}
\definecolor{darkpurple}{RGB}{128,0,128}

\usepackage[titletoc,toc,title]{appendix}

\newtheorem{theorem}{Theorem}[section]

\newtheorem{lemma}[theorem]{Lemma}
\newtheorem{cor}[theorem]{Corollary}

\theoremstyle{definition}

\newtheorem{rem}[theorem]{Remark}
\theoremstyle{definition}
\newtheorem{definition}[theorem]{Definition}
\newtheorem{example}[theorem]{Example}
\newtheorem{remark}[theorem]{Remark}

\def \cA {{\cal A}}
\def \cB {{\cal B}}

\def \cS {{\cal S}}

\def \Z {\mathbb Z}

\def \F {\mathbb F}

\begin{document}

\begin{frontmatter}
\title{On circular external difference families}
\author[1]{A.C.\ Burgess}
\ead{andrea.burgess@unb.ca}

\author[2]{F. Merola}
\ead{francesca.merola@uniroma3.it}

\author[3]{T.\ Traetta}
\ead{tommaso.traetta@unibs.it}

\affiliation[1]{organization={Department of Mathematics and Statistics, University of New Brunswick}, city={Saint John}, state={NB}, postcode={E2L 4L5}, country={Canada}}

\affiliation[2]{organization={Dipartimento di Matematica e Fisica, Università Roma Tre}, addressline={Largo S.L. Murialdo~1}, postcode={00142}, city={Roma}, country={Italy}}

\affiliation[3]{organization={DICATAM, Università degli Studi di Brescia}, addressline={Via Branze 43}, postcode={25123}, city={Brescia}, country={Italy}}



\begin{abstract}
Circular external difference families (CEDFs) are a recently-introduced variation of external difference families (EDFs) with applications to non-malleable threshold schemes:
a $(v,m,\ell,1)$-CEDF is an $m$-sequence $(A_0, \ldots, A_{m-1})$ of
$\ell$-subsets of an additive group $G$ of order $v$ such that
$G\setminus\{0\}$ equals the multiset
of all differences $a-a'$, with $(a,a')\in A_{i+1}\times A_{i}$ for some $i \in \mathbb{Z}_m$. When $G$ is the cyclic group, we speak of a cyclic CEDF.
 The existence of cyclic $(v,m,\ell,1)$-CEDFs is well understood when $m$ is even, while nonexistence is known when both $m$ and $\ell$ are odd. However, the case where $m$ is odd and $\ell$ is even has only been resolved in a few special cases.
 
 In this paper, we address this gap by constructing cyclic $(v,m,\ell,1)$-CEDFs  for any odd $m>1$ when $\ell=2$, and for any even $\ell \ge 2$ when $m=3$. 
Notably, the latter result relies on the existence of a suitable tiling of the multiplicative semigroup of $\mathbb{Z}_v\setminus\{0\}$. 
Moreover, noting that every $(v,3,\ell,1)$-CEDF produces a $(v,3,\ell,2)$-EDF, 
we completely solve the existence problem for 
$(3\ell^2+1,3,\ell,2)$-EDFs over an abelian group.

  Our approach is based on representing the blocks as arithmetic progressions and analyzing their step patterns. We present two different ways to construct cyclic $(v,m,2,1)$-CEDFs for every odd $m>1$; their step patterns show that the resulting CEDFs are inequivalent. Many additional inequivalent CEDFs are obtained by translating suitable subsets within the CEDF. 
\end{abstract}


\end{frontmatter}

\section{Introduction}
\emph{External Difference Families} (EDFs) have been intensively studied in the last 20 years, both for their combinatorial significance and for their applications to coding theory and cryptography \cite{CDFPW08, OKSS04}. In particular, there is a connection between  EDFs, some of their variations (see, for instance, \cite{JL}),  and  Algebraic Manipulation Detection Codes \cite{PS2}, which have applications, amongst others, to secret sharing schemes with special properties.
 
A new variation of EDFs,  \emph{Circular External Difference Families} (CEDFs), has been recently introduced in \cite{SV} as a tool to construct non-malleable threshold schemes. The definition is the following, denoting with $\Delta(A,B)$ 
the \emph{list of differences} between two subsets  $A$ and $B$  of a group $G$, that is, the multiset $\Delta(A,B) = \{a-b\mid a\in A,b\in B\}$:

\begin{definition}\label{def:CEDF}
Let $G$ be an additive group of order $v$. A sequence ${\cA}=(A_0,\dots, A_{m-1})$ of $m\ge 2$ disjoint $\ell$-sets is a $(v,m,\ell,\lambda)$-CEDF if the multiset union
\[\Delta(A_{1},A_0)\cup\Delta(A_2,A_1)\cup\dots\cup\Delta(A_{m-1},A_{m-2})\cup\Delta(A_{0},A_{m-1})\] 
is equal to
$\lambda (G\setminus \{0\}).$
If $G$ is cyclic, we speak of a cyclic CEDF.
\end{definition}
For instance, the sequence $\cA =
(\{0,1\},\{9,17\},\{3,6\},\{4,5\},\{16,18\})$ is a $(21,5,2,1)$-CEDF in the cyclic group $\Z_{21}$. 
Indeed,
\begin{align*}
   \Delta (\{9,17\},\{0,1\}) &= \{8,9,16,17\},\\
   \Delta (\{3,6\},\{9,17\}) &= \{7,10,15,18\},\\\
   \Delta (\{4,5\},\{3,6\}) &= \{1,2,19,20\},\\
   \Delta (\{16,18\},\{4,5\}) &= \{11,12,13,14\},\\
   \Delta (\{0,1\},\{16,18\}) &= \{3,4,5,6\},        
\end{align*}
so that the multiset union of the above lists of differences 
equals $\Z_{21}\setminus\{0\}.$
\begin{remark}
A more general notion (see \cite{SV}), not considered in this work, is that of a $c$-CEDF (for some positive integer $c$), say ${\cA}=(A_0,\dots, A_{m-1})$, where the property is that 
\[\textstyle{\bigcup_{i=0}^{m-1}\Delta(A_{i+c},A_{i})=\lambda (G\setminus \{0\})},\]
with the subscripts taken modulo $m$. 
When $c=1$, we obtain Definition~\ref{def:CEDF}. If $c$ and $m$ are relatively prime, the existence of a $c$-CEDF is equivalent to the existence of a CEDF.
\end{remark}

Note that CEDFs can be considered in the context of graph decompositions (\cite{HP, HJM}, see also \cite{BG}). If we denote by $\vec{C}_m[\ell]$ the lexicographic product of a directed $m$-cycle $\vec{C}_m$ with the empty graph on $\ell$ vertices, then a $(v,m,\ell,\lambda)$-CEDF is 
a vertex $G$-labeling $\Gamma$ of $\vec{C}_m[\ell]$ (that is, a digraph $\Gamma$ isomorphic to $\vec{C}_m[\ell]$, with $V(\Gamma)\subseteq G$) such that 
$\Delta \Gamma = \lambda(G\setminus\{0\})$. Note that, in this context, $\Delta \Gamma$ is the multiset of all differences $a-b$, provided that $(a,b)$ is an arc of $\Gamma$.
By using standard techniques, one can check that the existence of such a CEDF implies the existence of a $G$-regular decomposition of 
$\lambda K^*_v$ (i.e., the $\lambda$-fold symmetric complete digraph $K^*_v$) into copies of $\vec{C}_m[\ell]$.

A necessary condition for the existence of a $(v,m,\ell,\lambda)$-CEDF is that
$m\ell^2=\lambda(v-1)$ \cite{SV}, so that for $\lambda=1$ we have $v=m\ell^2+1$. CEDFs have been studied mainly when $\lambda=1$ and $G=\Z_v$ is the cyclic group of order $v$ (see \cite{PS,SV}); in particular, the results in \cite{PS} prove the existence of a cyclic $(v,m,\ell,1)$-CEDF  whenever $m$ is even. Also in $\cite{PS}$, Theorem~2.28 states the nonexistence of a cyclic $(v,m,\ell,1)$-CEDF when $m$ and $\ell$ are both odd (but see \cite{HJM} for examples in abelian, non-cyclic groups). The existence of a $(v,m,\ell,1)$-CEDF for $m$ odd and $\ell$ even is known only when $v=q$ is a prime power and $G=\F_q$, where  an extra condition must be satisfied
(\cite{SV}, see also Theorems 1.6 and 1.7 in \cite{PS}).
More precisely, we have the following.
\begin{theorem}[\cite{SV}] \label{PrimePowerCEDF}
Suppose that  $q=m\ell^2+1$ is a prime power and $\alpha$ is a primitive element of $\F_q$. Let $\beta=\alpha^\ell$ and let $H$ the subgroup of $\F_q^*$ of order $\ell m$ generated by $\beta$.  If 
$\{\beta-1,\beta^{m+1}-1,\dots,\beta^{(\ell-1)m+1}-1\}$
is a set of coset representatives of $H$ in $\F_q^*$,
then there exists a $(q,m,\ell,1)$-CEDF in $\F_q$.
\end{theorem}
In the case $\ell=2$, Theorem~\ref{PrimePowerCEDF} shows the existence of a $(4m+1,m,2,1)$-CEDF if $q=4m+1$ is a prime power and there exists a primitive element $\alpha$ of $\mathbb{F}_q$ such that $a^4-1$ is a quadratic non-residue in $\mathbb{F}_q^*$~\cite{SV}.  In~\cite{PS} (see also~\cite{WYF}), it is shown that whenever $q=4m+1$ is a prime power other than $5$, $9$ or $25$, there is a $(q,m,2,1)$-CEDF in $\mathbb{F}_q$.

The aim of this paper is to study cyclic CEDFs having $m$ odd and $\ell$ even, where $\lambda=1$. As mentioned above, the existence of such a CEDF is known only in some cases when $v=m\ell^2+1$ is a prime. 
In Section \ref{sec:l=2}, Theorem \ref{main:1} provides 
a $(v,m,\ell,1)$-CEDF for any odd $m\geq 3$ with $\ell=2$;
for $m=3$ and any even $\ell\geq 2$, existence is established in
Theorem \ref{main:2} (Section \ref{sec:m=3}). 
 In all cases, the proofs are constructive. The approach used to obtain these results is described in  Section~\ref{sec:approach}, 
where we formalize the notion of a CEDF whose sets are arithmetic progressions and explain how their ``steps'' play a central role in ensuring the CEDF properties. 

We now highlight a connection between classical and circular EDFs.
Recall that a $(v,m,\ell,\lambda)$-EDF over a group $G$ of order $v$
is a family ${\cA}=\{A_0,\dots,A_{m-1}\}$ of $m\ge2$ pairwise disjoint
$\ell$-subsets of $G$ such that 
$\bigcup_{i\neq j}\Delta(A_i,A_j)=\lambda(G\setminus\{0\})$.
Equivalent objects were introduced by Levenshtein \cite{Levenshtein} in 1971, under the name of
\emph{perfect and regular difference systems of sets}, for the construction of synchronization codes over cyclic groups.

In contrast to the case $m=2$, where cyclic EDFs are known to exist for
$\lambda=1,2$ (see, for instance, \cite{KPS}), only a finite number of EDFs
with $\lambda\leq 2$ and $m\geq 3$ seem to be currently known. In all these examples (see the table below), $m=3$ and
the underlying group is cyclic.

\begin{table}[h]
\centering
\renewcommand{\arraystretch}{1.2}
\begin{tabular}{|c|c|c|c|c|}
\hline
$(v,m,k,\lambda)$
& $(25,3,2,1)$
& $(97,3,4,1)$
& $(151,3,5,1)$
& $(13,3,2,2)$ \\
\hline
Reference
& \cite{GeMiaoWang2005}
& \cite{WangSu2010}
& \cite{CheeZhangZhang2011}
& \cite{FujiwaraTonchev2013} \\
\hline
\end{tabular}
\end{table}

Since the existence of a $(v,3,\ell,\lambda)$-CEDF immediately implies the
existence of a $(v,3,\ell,2\lambda)$-EDF, we obtain an infinite family
of $(v,3,\ell,2)$-EDFs;
indeed, Theorem~\ref{main:2} establishes the existence of a cyclic
$(3\ell^2+1,3,\ell,1)$-CEDF for every even $\ell\ge2$. Furthermore, in \cite[Theorem~5.1]{HJM}, the existence of $(3\ell^2+1,3,\ell,1)$-CEDF  over a suitable non-cyclic abelian group is shown whenever $\ell\equiv 3 \pmod{4}$, although the same construction given in this theorem also works when $\ell\equiv 1 \pmod{4}$.
We therefore obtain in Theorem~\ref{EDF} the existence  of a
$(3\ell^2+1,3,\ell,2)$-EDFs over some abelian group, for every $\ell\geq1$

In \cite{HJM}, the authors introduce the concept of equivalent CEDF, 
namely the case in which there is an affine transformation mapping one to another.
More precisely, two CEDFs over a group $G$, 
say $\cA=(A_{0},\ldots, A_{m-1})$ and $\cB=(B_{0},\ldots, B_{m-1})$, are said to be \emph{equivalent} if there is a triple $(\varphi, g, z)$, where $\varphi$ is an automorphism of $G$, $g\in G$, and $z\in \Z_m$ such that
$\varphi(A_i)+g = B_{i+z}$, for every $i\in\Z_m$. In the same paper \cite[Theorem 5.6]{HJM}, it is shown that there exist at least two inequivalent $(4\ell^2+1, 4, \ell, 1)$-CEDF for every non-prime $\ell>1$. To the best of our knowledge, this appears to be the only known result on inequivalent CEDFs.
Given the scarcity of results on inequivalent CEDFs, it is notable 
that the notion of ``pattern'' of an arithmetic CEDF (introduced in Section \ref{sec:approach}) can be used to easily show that the two CEDFs built in Theorems \ref{main:1} and \ref{main:0}  are inequivalent. Further results on inequivalent CEDFs will be presented in Section \ref{sec:inequivalent}. We conclude in Section \ref{sec:conclusions} with further questions and open problems.



\subsection{Preliminaries and notation.}
Here are some basic concepts and notation we will use throughout the paper.
Given two subsets $S_1, S_2$ of 
a group $G$, we define 
\[ 
\begin{aligned}
  S_1\pm S_2 =\{s_1\pm s_2\mid s_1\in S_1, s_2\in S_2\}.  
\end{aligned}
\] 
Clearly, $\Delta(S_1,S_2)=S_1-S_2$.  
If $S_i=\{s_i\}$, we simplify the notation by replacing $S_i$ with $s_i$, for $i=1,2$. 
Furthermore, 
letting  $x_1, x_2\in G$ and assuming that $G$ is abelian, it is clear that
\[\Delta(x_1+S_1, x_2+S_2) = (x_1-x_2) + \Delta (S_1, S_2).
\]
We also define the \emph{list of differences} of any subset $S$ of $G$ as the multiset $\Delta S = \{s-s'\mid s,s'\in S, s\neq s'\}$. Clearly, 
$\Delta (S+x) = \Delta S$.

Given the integers $x,y$ and $d\geq 1$, if $x\equiv y \pmod{d}$, 
we set
\[[x,y]_d = 
\begin{cases}
  \{x + id \mid 0\leq i\leq \frac{y-x}{d}\} & \text{if $x\leq y$},\\
  \varnothing & \text{if $x> y$}.
\end{cases}
\]
In the case $d=1$, we drop it from the notation, so when $x \leq y$, 
$[x,y]$ denotes the set of integers between $x$ and $y$, inclusive.

An \emph{arithmetic progression} of size $\ell\geq 2$ of $G$ 
(briefly, an $\ell$-AP)
is any $\ell$-set $A\subseteq G$ of the form
\[
A=[1, \ell]d + a = \{id+a\mid i\in [1,\ell]\},
\]
where $d, a\in G$ and $d\neq 0$. 
Recalling that the order $o(d)$ of $d$ equals the cardinality of the group $\langle d \rangle$ generated by $d$, we notice that $A$ is an $\ell$-set if and only if $\ell \leq o(d)$.
When $G=\Z_n$, this means that $\ell\leq o(d) = \frac{n}{\gcd(d,n)}$. 
Note also that $\Delta A = \Delta([1, \ell]d) = \Delta([1, \ell])\cdot d$.

Now, 
consider two $\ell$-AP of $G$, say 
$A=[1, \ell]d + a$ and $B=[1, \ell]f + b$.
If $\ell \leq  \min\{o(d), o(f)\} -2$, one can check that 
\[
A = B \
\iff 
(f,b) = (d,a)\;\;\text{or}\;\; (f,b) = (-d,(\ell+1)d+a).
\]
Indeed, $A=B$ implies that $\Delta A = \Delta B$. Since $\ell \leq \min\{o(d), o(f)\}-2$,
it is not difficult to check that $\pm d$ (resp. $\pm f$) are the only elements of the multiset $\Delta A$ (resp. $\Delta B$) with multiplicity $\ell-1$, hence $f=\pm d$, and the assertion easily follows. The converse is straightforward.

This justifies referring to $\delta(A)=\{d,-d\}$ as
the \emph{common difference} (or \emph{step}) of $A$; 
to simplify the notation, we will often write $\delta(A)=\pm d$.

\section{Our approach}
\label{sec:approach}
In \cite[Theorem 1.8 part (1)]{PS}, the authors build a cyclic $(m\ell^2+1,m,\ell,1)$-CEDF, say $\cS=(S_0, \ldots,$ $S_{m-1})$, for every even $m>0$. This result relies on a construction due to Snevily~\cite{Snevily} that builds $\alpha$-labelings of $G[n]$ (i.e. the lexicographic product of $G$ by the empty graph of order $n$) whenever $G$ has one.
As a consequence, the $\ell$-sets of $\cS$ 
are arithmetic progressions with steps  
$\pm 1$ and $\pm \ell$. 
Similarly, 
the $(3\ell^2+1,3,\ell,1)$-CEDF $\cA=(A_0, A_1, A_2)$ built in \cite[Theorem 5.1]{HJM} over $G=\Z_{\frac{3\ell^2+1}{2}}\times \Z_2$, for every  $\ell\equiv 3 \pmod{4}$, is another example where the sets forming the CEDF are arithmetic progressions; indeed,
\begin{align*}
  A_0 &= [0,\ell-1](1,0), \;\;\;\;\;\;
  A_1  = [0,\ell-1](z-1,1) + (-z-\ell,1),\\
  A_2 &= [0,\ell-1](z,1) + (-\ell,0),  
\end{align*}
with $z=\frac{3(\ell-1)^2}{4}$.
This motivates focusing on CEDFs with this property, which we call \emph{arithmetic};
in other words, 
a CEDF 
$\cA=(A_0, \ldots,A_{m-1})$ over $G$ is called arithmetic if each $A_i$ is an arithmetic progression of  $G$. In this case, any cyclic shift of
the sequence
$\delta(\cA)=(\delta(A_0), \ldots,$ $\delta(A_{m-1}))$
is called the \emph{pattern} of $\cA$, and the number $|supp(\delta(\cA))|$ of distinct entries of $\delta(\cA)$ is called the \emph{step-count}  of $\cA$. 
%
%
For example, the pattern of the arithmetic 
$(v,m,\ell,1)$-CEDF $\cS$ built in \cite{PS} 
for even $m$ is 
$\delta(\cS) = (\pm 1, \pm \ell, \ldots, 
\pm 1, \pm \ell)$, hence its step-count is $2$.  
Clearly, any CEDF is arithmetic when $\ell=2$. 

As shown below, CEDF with different step-counts are necessarily inequivalent.

\begin{lemma}\label{different_step_scounts}
Let $\cA=(A_{0},\ldots, A_{m-1})$ and $\cB=(B_{0},\ldots, B_{m-1})$ be arithmetic $(m\ell^2+1,m,\ell,\lambda)$-CEDFs over 
a  group $G$.  If $\cA$ and $\cB$ have different step-counts, then they are inequivalent.
\end{lemma}
\begin{proof}
We prove the contrapositive.  Let $(\varphi, g, z)$ be a triple such that $\varphi$ is an automorphism of $G$, $g\in G$, $z\in \Z_m$, and $\varphi(A_i)+g=B_{i+z}$ for each $i \in \mathbb{Z}_m$.  Given $i \in \mathbb{Z}_m$, if $A_i$ has common difference $\pm d$, then $B_{i+z}$ has common difference $\pm \varphi(d)$.  Hence if $\delta(A)=(\pm d_0, \pm d_1, \ldots, \pm d_{m-1})$ is the pattern of $\cA$, then the pattern of $\cB$ is \[\delta(\cB)=(\pm \varphi(d_{m-z}), \pm \varphi(d_{m-z+1}), \ldots, \pm \varphi(d_{m-1}), \pm \varphi(d_0), \ldots, \pm \varphi(d_{m-z-1})).\]  Since $\varphi$ is an automorphism, it is easy to see that $\delta(\cA)$ and $\delta(\cB)$ have the same number of distinct entries, i.e.\ $\cA$ and $\cB$ have the same step-count.
\end{proof}

\begin{example}
Consider the following $(21,5,2,1)$-CEDFs over $\mathbb{Z}_{21}$:
\[
\begin{aligned}
\mathcal{A}_1 &= \{ \{0,1\}, \{8,16\}, \{4,5\}, \{3,6\}, \{9,17\} \},  \\
\mathcal{A}_2 &= \{ \{0,1\}, \{16,18\}, \{4,5\}, \{3,6\}, \{9,17\} \}. \\
\end{aligned}
\]
Note that $\mathcal{A}_1$ and $\mathcal{A}_2$ have patterns $(\pm 1, \pm 8, \pm 1, \pm 3, \pm 8)$ and $(\pm 1, \pm 2, \pm 1, \pm 3, \pm 8)$, respectively, so that $\mathcal{A}_1$ has step-count $3$ and $\mathcal{A}_2$ has step-count $4$.  Thus, the CEDFs $\mathcal{A}_1$ and $\mathcal{A}_2$ are inequivalent.
\end{example}


Given two subsets $A$ and $B$ of an abelian group, it is known that $\Delta(A,B)$ (or, equivalently, $\Delta(B,A)$) has no repeated elements 
if and only if 
$\Delta A\,\cap\, \Delta B = \varnothing$. 
In particular, if $A$ and $B$ are arithmetic progressions with steps $\{\pm d_1\}$ and $\{\pm d_2\}$, respectively, this clearly implies that $d_1 \neq \pm d_2$. This condition turns out to be sufficient for $\Delta(A,B)$ to have no repeated elements only when $\ell=2$. 

In a CEDF $\cA=(A_0, \ldots, A_{m-1})$ with $\lambda=1$, each difference multiset
$\Delta(A_{i+1},A_i)$ cannot contain repeated elements; consequently 
$\delta(A_{i+1}) \neq  \delta(A_{i})$, 
for every
$i\in \Z_m$. It follows that the step-count of an arithmetic $(v,m,\ell,1)$-CEDF over $G$ is at least $3$ when $m$ is odd. Guided by these observations, and aiming to construct CEDF with patterns similar to the alternating one of
$\cS$ built in \cite{PS},
in Theorem \ref{main:1} we construct, for every odd $m\geq 3$, 
an (arithmetic) cyclic $(4m+1,m,2,1)$-CEDF with pattern
$(\pm1, \,\pm2, \,\pm1, \ldots, \,\pm2, \,\pm1, \,\pm3, \,\pm 2m-2)$
and step-count $4$. 
Since 
Theorem \ref{main:0} provides cyclic $(4m+1,m,2,1)$-CEDF with step-count $3$ (the least possible value), by Lemma \ref{different_step_scounts} these are inequivalent to those built in Theorem \ref{main:1}.

\section{The existence of cyclic $(4m+1,m,2,1)$-CEDFs}
\label{sec:l=2}
As mentioned above, the existence of a cyclic $(4m+1,m,2,1)$-CEDF with $m$ even has already been proven in \cite{PS}. The following theorem extends this result to all odd $m\geq3$, explicitly constructing a CEDF with step-count $4$ when $m> 3$.

\begin{theorem}\label{main:1}
There is a cyclic $(4m+1,m,2,1)$-CEDF for every odd $m \geq 3$.  
\end{theorem}
\begin{proof} Let $m=4u+1$ or $4u+3$ according to whether $m$ is congruent to 1 or 3 (mod 4), and for every $i \in \Z_m$, set $A_i= x_i + \{0,d_i\}$, where
\[
x_i=
\begin{cases}
  4m-2(i+1) & \text{if $i\in[1,2u-1]_2$},\\
  4(m-1)-2(i+1) & \text{if $i\in[2u+1,m-4]_2$},\\  
  2i & \text{if $i\in[0,m-3]_2$},\\
  2m-7  & \text{if $i=m-2$},\\
  2m-1   & \text{if $i=m-1$},   
\end{cases}
\]
and 
$(d_0, d_1, \ldots, d_{m-1}) = (1,2,1,\ldots,2,1,3,2m-2)$.
We claim that $\cA=(A_0,$ $A_1,$ $\dots,A_{m-1})$ is a $(4m+1,m,2,1)$-CEDF in $\Z_{4m+1}$. 

In the proof, we distinguish two cases.

\textbf{Case 1: $m\equiv 1 \pmod{4}$}

Let us first check that $\cA$ is a list of pairwise disjoint $2$-sets. Notice that
\[ 
\begin{aligned}
& \{x_0, x_2, \ldots, x_{m-3}\} = [0, 2m-6]_4, \\
& \{x_1, x_3, \ldots, x_{m-4}\} = [2m+2, 3m-7]_4\,\cup\,[3m+1, 4m-4]_4. 
\end{aligned}
\]
Therefore, letting $S:=\bigcup_{j=0}^{(m-3)/2} A_{2j}$ and $T:=\bigcup_{j=0}^{(m-5)/2} A_{2j+1}$,
we have that
\[
\begin{aligned}
  S &= \bigcup_{j=0}^{(m-3)/2} \{x_{2j}\} + \{0,1\} = 
  [0, 2m-6]_4 + \{0,1\} = [0, 2m-6]_4\,\cup\, [1, 2m-5]_4,\\
  T &= \bigcup_{j=0}^{(m-5)/2} \{x_{2j+1}\} + \{0,2\} 
    = \big([2m+2, 3m-7]_4\,\cup\,[3m+1, 4m-4]_4\big) + \{0,2\}  \\
   &= [2m+2, 3m-5]_2\,\cup\,[3m+1, 4m-2]_2. 
\end{aligned}  
\]
It is easy to see that $S$, $T$, $A_{m-2}= \{2m-7, 2m-4\}$ and $A_{m-1}=\{2m-1, 4m-3\}$ are pairwise disjoint.  Note that $|S| = m-1$ and $|T|= m-3$, so 
it follows that 
\[
\left|\bigcup_{i=0}^{m-1} A_i \right|  = 
\left| S\right| \,\cup\, 
\left| T\right| \,\cup\, 
\left|A_{m-2}\right| \,\cup\, 
\left|A_{m-1}\right|= (m-1) + (m-3) + 4 = 2m,
\]
and hence, the $A_i$s are pairwise disjoint. 

Letting $\Omega = \bigcup_{i=0}^{m-1}\Delta(A_{i+1}, A_{i})$, it remains to check that $\Omega=\Z_{4m+1}\setminus\{0\}$. We first notice that $\Delta(A_{i+1}, A_{i}) = x_{i+1}-x_{i} + \Delta(\{0, d_{i+1}\}, \{0, d_{i}\})$, where
\[
\begin{aligned}
\Delta(\{0, d_{i+1}\}, \{0, d_{i}\}) &=  \{0,d_{i+1}, -d_{i}, d_{i+1} -d_{i}\} \\
 &= \begin{cases} 
    [-1,2] & \text{if $i\in[0, m-5]_2$},\\
    [-2, 1] & \text{if $i\in[1, m-4]_2$},\\  
    \{-1,0,2,3\} & \text{if $i=m-3$},\\      
    \{-3, 0, 2m-2, 2m-5\} & \text{if $i=m-2$},\\     
    \{-2m+2,-2m+3,0,1\} & \text{if $i=m-1$},               
    \end{cases}
\end{aligned} 
\]
and 
\[
x_{i+1}-x_{i}=
\begin{cases}
  4(m-i-1) = -(4i+5) & \text{if $i\in[0,2u-2]_2$},\\
  4(m-i-2) = -(4i+9) & \text{if $i\in[2u, m-5]_2$},\\  
  4i+5 & \text{if $i\in[1,2u-1]_2$},\\
  4i+9 & \text{if $i\in[2u+1, m-4]_2$},\\ 
  -1=4m   & \text{if $i=m-3$},\\
  6   & \text{if $i=m-2$},\\
  2m+2 & \text{if $i=m-1$}.    
\end{cases}
\]
Letting 
$\Omega_0 = \bigcup_{j=0}^{(m-5)/2}{\Delta(A_{2j+1}, A_{2j})}$,
$\Omega_1 = \bigcup_{j=0}^{(m-5)/2}{\Delta(A_{2j+2}, A_{2j+1})}$ and
$\Omega_2 = \bigcup_{i=m-3}^{m-1} {\Delta(A_{i+1}, A_{i})}$,
it follows that 
\[
\begin{aligned}
 \Omega_0 &= 
    \bigcup_{j=0}^{(m-5)/2} \big( (x_{2j+1}-x_{2j}) + \Delta(\{0, d_{2j+1}\}, \{0, d_{2j}\}) \big)\\
   &= \bigcup_{j=0}^{(m-5)/2}  (x_{2j+1}-x_{2j}) + [-1,2]\\
   &= \big( [12,4m-8u-8]_8 \,\cup\, [4m-8u+4,4m-4]_8\big) +  [-1,2] \\
	 &= \big([12,2m-6]_8,\cup\, [2m+6,4m-4]_8\big)+[-1,2].
\end{aligned}
\]
and
\[
\begin{aligned}
 \Omega_1 &= 
     \bigcup_{j=0}^{(m-5)/2}  (x_{2j+2}-x_{2j+1}) + [-2,1]\\
   &= \big( [9,4m-8u-3]_8 \,\cup\, [4m-8u+9,4m-7]_8 \big) +  [-2,1]\\
	 &= \big([9,2m-1]_8 \,\cup\, [2m+11,4m-7]_8\big) + [-2,1].
\end{aligned}
\]
For the third set,
\[
\begin{aligned}
    \Omega_2 &=
      \big(4m+\{-1,0,2,3\}\big) \,\cup\,       
      \big(6+\{-3, 0, 2m-5, 2m-2\}\big) \,\cup \\ 
    &\;\;\;\;\,   
      \big(2m+2 + \{-2m+2,-2m+3,0,1\}\big) \\
    &= \pm \{1,2\} \,\cup\, \{3,6,2m+1,2m+4\} \,\cup\,   \{4,5,2m+2,2m+3\}  \\
    &= \pm \{1,2\} \,\cup\, [3,6] \,\cup\, [2m+1,2m+4].
\end{aligned}
\]
Therefore
\[
\begin{aligned}
    \Omega =&\; \Omega_0 \,\cup\, \Omega _1\,\cup\, \Omega_2 \\
=&\; \big([9,2m-1]_8+[-2,1]\big) \cup \big([12,2m-6]_8 + [-1,2]\big)\, \cup \\
&\; \big([2m+6,4m-4]_8 + [-1,2]\big) \cup \big([2m+11,4m-7]_8 + [-2,1] \big)\, \cup \\
&\; \pm \{1,2\} \,\cup\, [3,6] \,\cup\, [2m+1,2m+4] \\
=&\; \big([7,2m-3]_8 +[0,3]\big) \cup \big([11,2m-7]_8 + [0,3]\big) \,\cup \\
&\; \big([2m+5,4m-5]_8+[0,3]\big) \,\cup \big([2m+9,4m-9]_8+[0,3]\big) \,\cup \\
 &\; [1,2] \cup [4m-1,4m] \cup [3,6] \cup [2m+1,2m+4] \\
=&\; [7,2m] \cup [2m+5,4m-2] \cup [1,2] \cup [4m-1,4m] \cup [3,6] \,\cup\\
 &\; [2m+1,2m+4] \\
=&\; \mathbb{Z}_{4m+1} \setminus \{0\}.
\end{aligned}
\]
This completes the proof when $m\equiv 1\pmod{4}$.

\textbf{Case 2: $m\equiv 3 \pmod{4}$}

We reason as in the previous case, with minor modifications. The list $\cA$ consists of disjoint sets, since
\[ 
\begin{aligned}
& \{x_0, x_2, \ldots, x_{m-3}\} = [0, 2m-6]_4, \\
& \{x_1, x_3, \ldots, x_{m-4}\} = [2m+2, 3m-5]_4\,\cup\,[3m+3, 4m-4]_4. 
\end{aligned}
\]
Defining $S$ and $T$ as the previous case, we get
\[
\begin{aligned}
  S &= \bigcup_{j=0}^{(m-3)/2} \{x_{2j}\} + \{0,1\} = 
  [0, 2m-6]_4 + \{0,1\} = [0, 2m-6]_4\,\cup\, [1, 2m-5]_4,\\
  T &= \bigcup_{j=0}^{(m-5)/2} \{x_{2j+1}\} + \{0,2\} 
    = \big([2m+2, 3m-5]_4\,\cup\,[3m+3, 4m-4]_4\big) + \{0,2\}  \\ 
  &= [2m+2, 3m-3]_2\,\cup\,[3m+3, 4m-2]_2. 
 \end{aligned} 
\]
Once more,
considering that $A_{m-2} = \{2m-7, 2m-4\}$ and $A_{m-1} = \{2m-1, 4m-3\}$, 
it follows as before that the $A_i$s are pairwise disjoint also here.

Now define the sets $\Omega_0,\Omega_1,\Omega_2$ as above. The sets $\Omega_0$ and $\Omega_2$ are as in 
 the previous case:
\[
\begin{aligned}
 \Omega_0 &
     = \bigcup_{j=0}^{(m-5)/2}  (x_{2j+1}-x_{2j}) + [-1,2]\\
 &= \big( [12,4m-8u-8]_8 \,\cup\, [4m-8u+4,4m-4]_8\big) +  [-1,2] \\
 \end{aligned}
\]
and 
\[
    \Omega_2 = \pm \{1,2\} \,\cup\, [3,6] \,\cup\, [2m+1,2m+4],
\]
while for $\Omega_1$ we have
\[
\begin{aligned}
 \Omega_1 &= 
     \bigcup_{j=0}^{(m-5)/2}  (x_{2j+2}-x_{2j+1}) + [-2,1]\\
   &= \big([9,4m-8u-11]_8 \,\cup\, [4m-8u+1,4m-7]_8\big) +  [-2,1].
\end{aligned}
\]
Recalling that $u=(m-3)/4$, we see that in this case 
\[
\begin{aligned}
\Omega_0 
&= \big( [12,2m-2]_8 \,\cup\, [2m+10,4m-4]_8\big) + [-1,2] \\
&= \big([11,2m-3]_8 \, \cup \, [2m+9,4m-5]_8\big) + [0,3] 
\end{aligned}
\]
and
\[
\begin{aligned}
\Omega_1 &= \big([9,2m-5]_8 \,\cup\, [2m+7,4m-7]_8\big) + [-2,1] \\
&= \big([7,2m-7]_8 \, \cup \, [2m+5,4m-9]_8 \big) + [0,3].
\end{aligned}
\]
We show that $\Omega=\bigcup_{i=0}^2\Omega_i=\Z_{4m+1}\setminus\{0\}$ also in this case. Indeed,
\[
\begin{aligned}
    \Omega &= \Omega_0 \,\cup\, \Omega _1\,\cup\, \Omega_2\ \\
&= \big([7,2m-7]_8 + [0,3]\big) \cup \big([11,2m-3]_8 + [0,3]\big) \cup \\
&\;\;\;\;\;\; \big([2m+5,4m-9]_8+[0,3]\big) \cup \big([2m+9,4m-5]_8+[0,3]\big) \cup \\
&\;\;\;\;\;\; \pm \{1,2\} \,\cup\, [3,6] \,\cup\, [2m+1,2m+4] \\
&= [7,2m] \cup [2m+5,4m-2] \cup [1,2] \cup [4m-1,4m] \cup [3,6] \cup \\
&\;\;\;\;\;\; [2m+1,2m+4] \\
&= \mathbb{Z}_{4m+1} \setminus \{0\}.
\end{aligned}
\]
This completes the proof.
\end{proof}
We note that Theorem~\ref{main:1} constructs a CEDF with pattern 
$(\pm 1, \pm 2, \ldots,$ 
$\pm 1,\pm 2, \pm 1, \pm 3, \pm (2m-2))$; hence its step-count is $4$ when $m > 3$.  
\begin{cor}
For all $m \geq 5$, there is a cyclic $(4m+1,m,2,1)$-CEDF with step-count $4$.
\end{cor}

\begin{example}
Using the construction above with $m=7$ we obtain the following 
CEDF in $\Z_{29}$:
\[\cA=(\{0,1\},\{24,26\},\{4,5\},\{16,18\},\{8,9\},\{7,10\},\{13,25\}).\]
With $m=9$ we have the following CEDF in $\Z_{37}$:
\[\cB=(\{0,1\},\{32,34\},\{4,5\},\{28,30\},\{8,9\},\{20,22\},\{12,13\},\{11,14\},\{17,33\}).\]
Since their patterns take the form
$(\pm 1,\pm 2, \pm1,\ldots,\pm 2,\pm 1,\pm 3,\pm (2m-2))$, they both have step-count 4.
\end{example}

\begin{rem} 
For every odd $m\geq 3$, there is a a close relationship between the
$(4m-3,m-1,2,1)$-CEDF $\cS$ constructed in \cite{PS} and the
$(4m+1,m,2,1)$-CEDF $\cB$ constructed in Theorem~\ref{main:1}.
To illustrate this connection, consider the case $m=9$. 
The $(37,9,2,1)$-CEDF $\cB$
is given in the previous example, whereas $(33,8,2,1)$-CEDF $\cS$ is
\[
\cS=
(\{0,1\},\{30,32\},\{4,5\},\{26,28\},\{8,9\},
 \{18,20\},\{12,13\},\{14,16\}).
\]
Denote by $B_i$ the $i$-th block of $\cB$ and by $S_i$ the $i$-th
block of $\cS$. Then
\[
B_i=
\begin{cases}
S_i,   & \text{if $i$ is odd},\\[1mm]
S_i+2, & \text{if $i$ is even},
\end{cases}
\qquad 1\leq i\leq 7.
\]
Consequently, for every $1\leq i\leq 6$, the list of differences in
$\Z$ from $B_{i+1}$ to $B_{i}$ is obtained from the corresponding list
of differences from $S_{i+1}$ to $S_{i}$ by adding $2$ when $i$ is odd
and subtracting $2$ when $i$ is even. More precisely,
\[
\Delta(B_{i+1},B_i)=
\begin{cases}
\Delta(S_{i+1},S_i)+2, & \text{if $i$ is odd},\\[1mm]
\Delta(S_{i+1},S_i)-2, & \text{if $i$ is even}.
\end{cases}
\]
Setting
$\Delta_{\cB}:=
\bigcup_{i=1}^{6}\Delta(B_{i+1},B_i)$,
and working modulo $37$, the only nonzero residues that do not occur in
$\Delta_{\cB}$ are
\[
-1,\,-2,\,2,\,1,\,
6,\,3,\,22,\,19,\,
-17,\,-33,\,-16,\,-32.
\]
These missing residues are then precisely covered by 
\[
\Delta(B_8,B_7)=\{-1,-2,2,1\},
\quad
\Delta(B_9,B_8)=\{6,3,22,19\},\quad\text{and}
\]
\[
\Delta(B_1,B_9)=\{-17,-33,-16,-32\}.
\]
\end{rem}

Theorem~\ref{main:1}, combined with Theorem 1.8 of~\cite{PS}, completely settles the existence of cyclic $(v,m,\ell,1)$-CEDFs with $\ell=2$.
\begin{theorem}
Let $v,m>1$ be integers.  There exists a cyclic $(v,m,2,1)$-CEDF if and only if $v=4m+1$.
\end{theorem}

\section{The existence of cyclic $(3\ell^2+1,3,\ell,1)$-CEDFs}
\label{sec:m=3}
In this section, we prove the second main result of this paper, Theorem~\ref{main:2}, where we build a $(3\ell^2+1,3,\ell,1)$-CEDF over $\Z_{3\ell^2+1}$ whenever the trivial necessary condition for its existence holds, that is, for every even $\ell>0$.
\begin{theorem}\label{main:2}
There exists an arithmetic $(3\ell^2+1,3,\ell,1)$-CEDF in $\Z_{3\ell^2+1}$ 
for every even $\ell\geq 2$. 
\end{theorem}
This result is proven at the end of this section where we explicitly provide the three sets $A_1, A_2, A_3$ constituting the desired CEDF and show that
\[
\textstyle{\bigcup_{i\in\Z_3}} \Delta (A_{i+1}, A_i) = R\cdot\{1,d,d^2\},
\]
where $R=\big([1,2k]d + [0,2k-1]\big)$ and $d=6k^2-3k$.
We then use Lemmas \ref{lemma:21}  and \ref{factorization} to show that $R\cdot\{1,d,d^2\} = \Z_v\setminus\{0\}$.

Throughout this section, we assume that 
\[
\text{$\ell=2k$,\;\;\; $d=6k^2-3k$\;\;\; and\;\;\;  $v=3\ell^2+1 = 12k^2+1$,}
\] 
for some positive integer $k$.
Lemma \ref{lemma:21} determines all possible integer solutions $(\alpha,\beta)$ to the congruence equation
\begin{equation}\label{cong_eq}
  dX+Y \equiv 0 \pmod{v}
\end{equation}
whenever $(\alpha,\beta)$ is constrained to belong to some specific subsets of $\Z^2$.
Denote by $\Sigma=\{(\alpha,-d\alpha+hv)\mid \alpha,h\in \Z\}$ the set of integral solutions of \eqref{cong_eq}.
Note that any integer $\alpha$ can be uniquely expressed  in the following form 
\[\alpha= (2k-1)a + 2b + \epsilon,
\]
for some $a\in \Z$, $b\in[0,k-1], \epsilon\in\{0,1\}$, and $(b,\epsilon)\neq (k-1,1)$;
also, set
\[
  \psi(\alpha) = -(2k+1)a + (6k+1)(b+\epsilon) + d\epsilon.
\]
Considering that $v=2d+6k+1$ and $v(k-1)=d(2k-1)-(2k+1)$, one can easily check that 
$\psi(\alpha) \equiv -d\alpha \pmod{v}$, for every $\alpha\in \Z$; indeed,
\[\begin{aligned}
-d\alpha &= -(2k-1)ad -2bd - \epsilon d  \\
&\equiv-(2k+1)a + (6k+1)b +  (d+6k+1)\epsilon \pmod{v}\\
&= \psi(\alpha).
\end{aligned}
\]
Therefore, the set $\Sigma$ of integral solutions to \eqref{cong_eq} can be written as follows:
\[\Sigma=\{(\alpha, \psi(\alpha)+hv)\mid \alpha,h\in \Z\}.\]
\begin{lemma} \label{lemma:21} 
The only integral solutions to \eqref{cong_eq} in the set $[1-2k, 4k]\times [0,4k-1]$ are $(1-2k,2k+1),(0,0)$, and $(4k, 2k-1)$.
\end{lemma}
\begin{proof} Let $(\alpha,\beta)\in [1-2k, 4k]\times [0,4k-1]$ be an integral solution of \eqref{cong_eq}. 
Since $\alpha\in [1-2k, 4k]$, it can be uniquely expressed in the form 
\[\alpha = (2k-1)a + 2b + \epsilon,\] 
for some $a\in [-1,2]$, $b\in[0,k-1]$ and $\epsilon\in\{0,1\}$ such that
\begin{equation}\label{eq:22}
\begin{aligned}
(b,\epsilon)\neq (k-1,1)
\end{aligned}
\end{equation}
Letting $\psi(\alpha) = -(2k+1)a + (6k+1)(b+\epsilon) + d\epsilon$, we start by showing 
that $\beta=\psi(\alpha)$; since 
\begin{equation}\label{eq:23}
\text{$0\leq \beta < 4k$\;\;\; and\;\;\; $\beta\equiv \psi(\alpha) \pmod{v}$,}
\end{equation}
it is enough to show that 
$0\leq \psi(\alpha) < v$.
Indeed, if $\psi(\alpha)< 0$, then $a\in\{1,2\}$ and $b=\epsilon=0$, that is,
$\alpha\in \{2k-1, 4k-2\}$. Since 
\begin{equation}\label{eq:21}
 \text{$\psi(2k-1) = -(2k+1)$,\;\;\; $\psi (4k-2) = -(4k+2)$,}
\end{equation}
and in view of \eqref{eq:23}, it follows that
$4k> \beta\geq v+\psi(\alpha)$, that is, $\psi(\alpha)<4k-v < -(4k+2)$, thus contradicting
\eqref{eq:21}.  Similarly, if $\psi(\alpha)\geq v$, then 
$a=-1$ and $(b,\epsilon)=(k-1,1)$, thus contradicting \eqref{eq:22}.

Therefore, $0\leq \beta=\psi(\alpha) < 4k$. Note that if $\epsilon=1$, we would have
\[ 
  \psi(\alpha) = d-(2k+1)a + (6k+1)(b+1) \geq d-2(2k+1) + (6k+1) =6k^2-k-1\geq 4k.
\] 										
Therefore, $\epsilon = 0$ and
$\beta =  \psi(\alpha) = -(2k+1)a + (6k+1)b$. 
One can easily check that
$0\leq \psi(\alpha)< 4k$ if and only if 
$(a,b)\in \{(-1,0),(0,0), (2,1)\}$. This is equivalent to saying that
$(\alpha,\beta) \in \{(1-2k,2k+1),(0,0), (4k, 2k-1)\}$, and this completes the proof.
\end{proof}

The following result provides a ``tiling'' of the multiplicative monoid of $\Z_v$.

\begin{lemma} \label{factorization}
Let $d=6k^2-3k$ and $R=[1,2k]d+[0,2k-1]$. Then, in $\Z_{12k^2+1}$, we have that
$d^3=1$, $d^2=-(d+1)$ and $R\cdot\{1, d, d^2\} = \Z_{12k^2+1}\setminus\{0\}$.
\end{lemma}
\begin{proof} Let $v=12k^2+1$. In $\Z_v$, it is easy to check that $d^2=-(d+1)$ and $d^3=1$, which implies that $d$ is invertibile in $\Z_v$. Therefore, to prove that
$R\cdot\{1, d, d^2\} = \Z_{12k^2+1}\setminus\{0\}$ it is enough to show that $0\not\in R$  and
$R\,\cap\, (R\cdot d) = \varnothing$.

Clearly, $0\not\in R$; otherwise there exist $\alpha\in[1,2k]$ and
$\beta\in [0,2k-1]$ such that $d\alpha + \beta\equiv 0\pmod{v}$, thus contradicting 
Lemma \ref{lemma:21}. Now, assume for a contradiction that 
$R\,\cap\, (R\cdot d) \neq \varnothing$, which is equivalent to saying that there exist $\alpha, \alpha'\in[1,2k]$ and $\beta,\beta'\in [0,2k-1]$ such that
\[
  d\alpha + \beta 
  \equiv (d\alpha' + \beta')d 
  \equiv d^2\alpha' + d \beta' 
  \equiv -(d+1)\alpha' + d \beta'\pmod{v}.
\]
Hence 
$d(\alpha + \alpha'  - \beta')  + \alpha' + \beta\equiv 0\pmod{v}$, 
that is, $(\alpha + \alpha'  - \beta', \alpha' + \beta) \in [3-2k,4k]\times[1,4k-1]$
is an integral solution to $dX+Y\equiv 0\pmod{v}$. By Lemma \ref{lemma:21}, 
$(\alpha + \alpha'  - \beta', \alpha' + \beta)= (4k,2k-1)$. Therefore,
\[  \alpha - \beta'  - \beta = 
(\alpha + \alpha'  - \beta') - (\alpha' + \beta) = 2k+1,
\]
but this is impossible since $\alpha - \beta'  - \beta \leq 2k$.
\end{proof}

We are now ready to prove the main result of this section.

\begin{proof}[Proof of Theorem \ref{main:2}]  Let $\ell=2k$, $d= 6k^2-3k$ and $v=3\ell^2+1=12k^2+1$, for some $k\geq 1$. Considering that the case $\ell=2$ is solved in Theorem \ref{main:1}, we can assume that $k\geq2$.
We claim that $\cA=(A_0,A_1,A_{2})$, where
\[
\begin{aligned}
	A_0 &= [1, 2k],\;\;\;
	A_1 &= [1,2k]d + 2k, \;\;\;\text{and}\;\;\;
	A_2 &= [1,2k]d^2 + 6k^2 + k + 1, 
\end{aligned}
\]
is an arithmetic $(v,3,\ell,1)$-CEDF in $\Z_{v}$. Clearly, each $A_i$ is a $2k$-AP.
Considering that $6k^2 - k + 1= 2kd$ and $-(6k^2 + k + 1) = 2kd^2$ in $\Z_v$ and recalling that $d^3=1$ so that $[1,2k]=[1,2k]d^3$, 
it follows that
\[
\begin{aligned}
	\Delta (A_1, A_0) &= [1,2k]d -[1, 2k] + 2k =  [1,2k]d + [0,2k-1],	\\
	\Delta (A_2, A_1) &= [1,2k]d^2 - [1,2k]d + 6k^2 - k + 1  
	                  = [1,2k]d^2 - [1,2k]d + 2kd  \\
	                  &=  ([1,2k]d + [0,2k-1])d, \\
	\Delta (A_0, A_2) &= [1, 2k]d^3 - [1,2k]d^2 - (6k^2 + k + 1) 
	                  = [1, 2k]d^3 - [1,2k]d^2 + 2kd^2 \\
	                  &= ([1,2k]d + [0,2k-1])d^2.
\end{aligned}
\]
Therefore, by Lemma \ref{factorization}, we have that
\[
\textstyle{\bigcup_{i\in\Z_3}} \Delta (A_{i+1}, A_i) = \big([1,2k]d + [0,2k-1]\big)\cdot\{1,d,d^2\}
= \Z_v\setminus\{0\}.\]
Since $\cA$ consists of three sets and each $\Delta (A_{i+1}, A_i)$ 
does not contain $0$, we have that $A_0, A_1$ and $A_2$ are pairwise disjoint, and this completes the proof.
\end{proof}

\begin{example}
Letting $\ell=4$, we have that $k=2, d=18$ and $v=49$. By Theorem \ref{main:2}, we obtain a $(49,3,4,1)$-CEDF in $\mathbb{Z}_{49}$, say $(A_0, A_1, A_2)$,  where:
\begin{align*}
  A_0 &= {\{1,2,3,4\},} \\ 
  A_1 &= {\{22,40,9,27\},} \\ 
  A_2 &= {\{8,38,19,0\}.}  
\end{align*} 
Letting $\ell=6$, we have that $k=3, d=45$ and $v=109$. By Theorem \ref{main:2}, we obtain a $(109,3,6,1)$-CEDF in $\mathbb{Z}_{109}$, say $(A_0, A_1, A_2)$, where:
\begin{align*}
  A_0 &=  {\{1,2,3,4,5,6\},} \\ 
  A_1 &=  {\{51,96,32,77,13,58\},} \\ 
  A_2 &=  {\{12,75,29,92,46,0\}.}  
\end{align*} 
Letting $\ell=8$, we have that $k=4, d=84$ and $v=193$. By Theorem \ref{main:2}, we obtain a $(193,3,8,1)$-CEDF in $\mathbb{Z}_{193}$, say $(A_0, A_1, A_2)$, where:
\begin{align*}
  A_0 &= {\{1,2,3,4,5,6,7,8\},} \\ 
  A_1 &= {\{92,176,67,151,42,126,17,101\},} \\ 
  A_2 &= {\{16,124,39,147,62,170,85,0\}.}  
\end{align*} 
\end{example}

Theorem~\ref{main:2}, combined with Theorem 1.8 of~\cite{PS}, completely settles the existence of cyclic $(v,m,\ell,1)$-CEDFs with $m=3$.
\begin{theorem}
Let $v,\ell\geq 1$ be integers.  There exists a cyclic $(v,3,\ell,1)$-CEDF if and only if $v=3\ell^2+1$ and $\ell$ is even.
\end{theorem}

As already noted in the introduction, the existence of a $(v,3,\ell,\lambda)$-CEDF implies the existence of a $(v,3,\ell,2\lambda)$-EDF.
Therefore, by Theorem~\ref{main:2} and recalling that the construction in Theorem~5.1 of \cite{HJM} proves the existence of a $(v,3,\ell,1)$-CEDF over a suitable non-cyclic abelian group for every odd $\ell\geq 1$, we obtain the following result.
\begin{theorem}\label{EDF}
There exists a $(3\ell^2+1,3,\ell,2)$-EDF over an abelian group for every $\ell\geq 1$.
\end{theorem}

\section{Inequivalent CEDFs}\label{sec:inequivalent}
In this section, we present Theorem~\ref{main:0}, which provides an alternative construction for a cyclic $(4m+1,m,2,1)$-CEDF, this time with step-count 3.
Recall that two CEDFs over a group $G$, 
say $\cA=(A_{0},\ldots, A_{m-1})$ and $\cB=(B_{0},\ldots, B_{m-1})$, are said to be equivalent if there is a triple $(\varphi, g, z)$, where $\varphi$ is an automorphism of $G$, $g\in G$, and $z\in \Z_m$ such that
$\varphi(A_i)+g = B_{i+z}$, for every $i\in\Z_m$. Since two CEDFs with distinct step-counts are necessarily inequivalent (Lemma \ref{different_step_scounts}), it follows that these CEDFs are inequivalent from the CEDFs with step-count $4$ constructed in Theorem~\ref{main:1}.  We will then describe a method that, by a slight modification of a given arithmetic $(v,m,\ell,\lambda)$-CEDF allows us to possibly obtain many inequivalent CEDFs with the same parameters.

\begin{theorem}\label{main:0}
There is a cyclic $(4m+1,m,2,1)$-CEDF with step-count $3$, for every odd $m>1$.
\end{theorem}

\begin{proof} We will build CEDFs with pattern 
$(\pm1, \pm (2m-2), \pm 1, \ldots, \pm(2m-2), \pm 1, \pm 3, \pm (2m-2))$, and so they will have step-count $3$. For $m=3,5,7$, the desired CEDFs are given as follows: 
\begin{align}
&\{ \{0,1\}, \{12,2\}, \{5,9\} \}, 
\label{small_cases_1} \\[2pt]
&\{ \{0,1\}, \{8,16\}, \{4,5\}, \{3,6\}, \{9,17\} \}, \label{small_cases_2} \\[2pt]
&\{ \{0,1\}, \{14,26\}, \{4,5\}, \{16,28\}, \{8,9\}, \{7,10\}, \{13,25\} \}.
\label{small_cases_3} 
\end{align}
Now, let $m=8b+2\epsilon+1$, with $\epsilon\in[0,3]$, and assume that
$m\geq 9$, that is, $b\geq 1$. 
For every $i \in \Z_m$, set $A_i= x_i + \{0,d_i\}$, where
$(d_0, d_1, \ldots, d_{m-1}) = (1,2m-2,1,\ldots,2m-2,1,3,2m-2)$
and
\[
x_i=
\begin{cases}
  2i & \text{if $i\in[0,m-3]_2$},\\
  2m-2i+6    & \text{if $i\in[3,4b-1]_4$},  \\
  2m-2i-10   & \text{if $i\in[1,4b-7]_4$},  \\  
  2m-8b-4\lceil \frac{\epsilon}{4}\rceil
  & \text{if $i=4b-3$},  \\
  4m-2i-10 
  & \text{if $i\in[4b+1, m-4]_2$},\\
   2m-7  & \text{if $i=m-2$},\\
   2m-1   & \text{if $i=m-1$};  \\
\end{cases}
\]
also, when $\epsilon=3$ (i.e.\ $m \equiv 7 \pmod{8}$), we replace two of these values by setting 
\[(x_{4b+1}, x_{4b+3}) = (m+7,\, 3m-5) = (8b+14, 24b+16).\]
Note that when $i=4b-3$ and $\epsilon \neq 0$, $x_i = 2m-2i-10$; if $\epsilon=0$, $x_{i} = 2m-2i-6$. 
We claim that $\cA=(A_0,A_1,\dots,A_{m-1})$ is a $(4m+1,m,2,1)$-CEDF in $\Z_{4m+1}$; this is verified in \ref{Appendix}.  
\end{proof}

Readers who prefer to skip the technical details of the proof of Theorem~\ref{main:0} may find it helpful to instead consider the underlying idea of the construction, which is inspired by the examples in
\eqref{small_cases_1} and \eqref{small_cases_2}, whose patterns are
\[
(1,3,2m-2)
\quad\text{and}\quad
(1,2m-2,1,3,2m-2),
\]
respectively. The general construction preserves the same local structure.
More precisely, the sequence 
$\mathcal A=(A_0,A_1,\ldots,A_{m-1})$ (representing the CEDF in Theorem~\ref{main:0}) 
can be naturally decomposed into the terminal section
$S=(A_{m-3}$, $A_{m-2},A_{m-1},A_0)$ whose pattern is $(1,3, 2m-2,1)$,
and the pairs of consecutive sets $(A_i,A_{i+1})$, for $i\in [0,m-4]$, whose patterns are either $(1, 2m-2)$ or $(2m-2,1)$.
The role of the subsequence $S$ is to generate a small exceptional set of differences. Indeed, since 
\begin{align*} 
A_{m-3}&=\{2m-6,2m-5\},\quad 
A_{m-2}=\{2m-7,2m-4\},\quad \\
A_{m-1}&=\{2m-1, 4m-3\},\quad 
A_{0}=\{0,1\}, 
\end{align*}
it is not difficult to check that 
\[\Delta S = [1,6]\,\cup\,[2m+1,2m+4]\,\cup\,\{4m-1,4m\}.\] 
Each of the remaining pairs $(A_i, A_{i+1})$
contributes four differences arranged into two pairs of the form 
$\{x, x+2m-2\}$, thus allowing us to partition (as differences) the remaining integers in $\Z_{4m+1}\setminus\{0\}$ not covered by $\Delta S$, that is, the set 
\[\Z_{4m+1}\setminus (\Delta S \,\cup\,\{0\})= [7,2m]\,\cup\,[2m+5, 4m-2] = [7,2m]+\{0, 2m-2\}.\] 
Indeed, since 
\[A_i = 
\begin{cases} \{a,a+1\} & \text{if $i\in[0,m-4]$ is even},\\ 
\{b,b+1\} & \text{if $i\in[0,m-4]$ is odd}, \end{cases}
\] 
each $\Delta (A_i, A_{i+1})$ takes the following form \[ \Delta (A_i, A_{i+1})= 
\begin{cases} 
  (b-a) + \{-1,0, 2m-3, 2m-1\} & \text{if $i$ is even},\\ 
  (a-b) + \{-2m+2, -2m+3, 0, 1\} & \text{if $i$ is odd}. 
\end{cases} 
\]

\begin{example}
Using the construction of Theorem~\ref{main:0} with $m=9$, we obtain the following CEDF in $\mathbb{Z}_{37}$:
\[
\mathcal{A}=(\{0,1\}, \{10,26\}, \{4,5\}, \{18,34\}, \{8,9\}, \{16,32\}, \{12,13\}, \{11,14\}, \{17,33\}).
\]
This $(37,9,2,1)$-CEDF has pattern 
$(\pm1,\pm16, \ldots,\pm1,\pm16,\pm1,\pm3,\pm16)$.

Taking $m=15$, we obtain the following $(61,15,2,1)$-CEDF in $\mathbb{Z}_{61}$, which has pattern 
$(\pm 1,\pm 28,\ldots, \pm 1,\pm 28, \pm1, \pm 3,\pm 28)$:
\[
\begin{array}{l}
\mathcal{A}=(\{0,1\}, \{18,46\}, \{4,5\}, \{30,58\}, \{8,9\}, \{22,50\}, \{12,13\}, \{40,7\}, \{16,17\}, \\
\hspace*{1.5in} \{32,60\}, \{20,21\}, \{28,56\}, \{24,25\}, \{23,26\}, \{29,57\}).
\end{array}
\]
\end{example}

We now present a method for modifying specific sets within a CEDF while preserving the overall list of differences. Therefore, if the resulting sequence consists of pairwise disjoint sets, we obtain a new CEDF. 
Before proceeding, we give an example to illustrate the idea of the method.  Consider the cyclic $(61,15,2,1)$-CEDF 
$\mathcal{A}=(A_0,A_1, \ldots, A_{14})$ from the previous example.
Note that $A_2-A_1 = \{19,20,47,48\}$, while $A_3-A_2=\{25,26,53,54\} = A_2-A_1 + 6$.  If we add 6 to each element of $A_2=\{4,5\}$ to obtain $A_2'=\{10,11\}$ we see that $A_2'-A_1 = A_2-A_1+6 = A_3-A_2$ and $A_3-A_2' = A_3-A_2-6=A_2-A_1$.  Consequently, $\Delta(A_2,A_1) \cup \Delta(A_3,A_2) = \Delta(A_2',A_1) \cup \Delta(A_3,A_2')$.  Since $A_2'$ is disjoint from each $A_i$, $i \in \{0, \ldots, 14\} \setminus \{2\}$, replacing $A_2$ with $A_2'$ yields a $(61,15,2,1)$-CEDF $\mathcal{A}'$, given below:
\[
\begin{array}{l}
\mathcal{A}'=(\{0, 1\}, \{18, 46\}, \{10, 11\}, \{30, 58\}, \{8, 9\}, \{22, 50\}, \{12, 13\}, \{40, 7\}, \\
\hspace*{1cm}\{16, 17\}, \{32, 60\}, \{20, 21\}, \{28, 56\}, \{24, 25\}, \{23, 26\}, \{29, 57\}).
\end{array}
\]
Similarly, replacing $A_4=\{8,9\}$ in $\mathcal{A}$ with $A_4'=A_4+2=\{10,11\}$ also preserves the overall list of differences, and since $A_4'$ is disjoint from $A_i$ for $i \in \{0, \ldots, 14\} \setminus \{4\}$, replacing $A_4$ in $\mathcal{A}$ with $A_4'$ also yields a $(61,15,2,1)$-CEDF $\mathcal{A}''$:
\[
\begin{array}{l}
\mathcal{A}''=(\{0, 1\}, \{18, 46\}, \{4, 5\}, \{30, 58\}, \{10,11\}, \{22, 50\}, \{12, 13\}, \{40, 7\}, \\
\hspace*{1cm}\{16, 17\}, \{32, 60\}, \{20, 21\}, \{28, 56\}, \{24, 25\}, \{23, 26\}, \{29, 57\}).
\end{array}
\]
Note, however, that $(A_0,A_1,A_2',A_3,A_4',A_5,\ldots,A_{14})$ is not a CEDF; while the overall list of differences is preserved, $A_2'$ and $A_4'$ are not disjoint.  

We now describe the method more precisely.
Let $\mathcal{A}=(A_0, \ldots, A_{m-1})$ be an arithmetic $(m\ell^2+1, m, \ell, \lambda)$-CEDF over an abelian group $G$, where each $A_i = [0, \ell-1]d_i+x_i$. 
Furthermore, 
assume there exist $0\leq i\leq j \leq m-1$, with $(i,j)\neq (0,m-1)$, such that
$\pm \{d_{i-1}, d_i\} = \pm \{d_{j}, d_{j+1}\}$, and  set
\[
 z= x_{j+1}-x_{j} -(x_{i}-x_{i-1})+
 \begin{cases}
   0                     &  \text{if $(d_{i-1}, d_i) = (d_{j}, d_{j+1})$},\\
   (d_{j+1}-d_j)(\ell-1) &  \text{if $(d_{i-1}, d_i) = (d_{j+1}, d_{j})$}.
 \end{cases} 
\]
We say that the sequence $\mathcal{A}' = (A'_0, \ldots, A'_{m-1})$ where
$A_h' = A_h +z$ if $h\in [i,j]$, and $A_h' = A_h$ otherwise, 
is \emph{$(i,j)$-associated} to $\cA$ (or simply \emph{$i$-associated}, if $i=j$). One can check that
\[
 \Delta (A_{j+1}, A_j)\,\cup\, \Delta (A_{i}, A_{i-1})
 =
 \Delta (A_{j+1}, A_j+z)\,\cup\, \Delta (A_{i}+z, A_{i-1}),
\]
thus implying that $\bigcup_{h\in \Z_m} \Delta(A_{h+1},A_h) = 
\bigcup_{h\in \Z_m} \Delta(A'_{h+1},A'_h)$. We then have the following result.

\begin{lemma} \label{lem:associatedCEDF}
  Let $\mathcal{A}=(A_0, \ldots, A_{m-1})$ be an arithmetic CEDF over an abelian group $G$ such that
\[\delta(A_{j-1})\,\cup\, \delta(A_{j}) = \delta(A_{i})\,\cup\, \delta(A_{i+1}),\] 
with $0\leq i\leq j \leq m-1$ and $(i,j)\neq (0,m-1)$. Also, let $\cA'$ be
the sequence $(i,j)$-associated to $\cA$. Then, $\cA'$ is a CEDF if and only if its sets are pairwise disjoint.
\end{lemma}

In Theorem~\ref{th:inequivalent}, we will construct, for each odd integer $m \geq 33$, $4\lfloor m/24 \rfloor-3$ pairwise inequivalent $(4m+1,m,2,1)$-CEDFs with the same pattern by finding CEDFs which are $i$-associated to the CEDF constructed in Theorem~\ref{main:0} for appropriate values of $i$.  To prove that these CEDFs are indeed inequivalent, we will make use of the following lemma.
\begin{lemma}\label{lem:ineq:1}
  Let $\cA$ and $\cB$ be two cyclic arithmetic $(m\ell^2+1,m,\ell,\lambda)$-CEDF 
  with pattern $(\pm 1, \pm d,\ldots, \pm 1,  \pm3, \pm d)$. 
  If $\cA\neq \pm \cB$ and $3(d-1)\neq0$ in $\Z_{v}$, then $\cA$ and $\cB$ are 
  inequivalent.
\end{lemma}
\begin{proof}
  Let $\cA=(A_{0},\ldots, A_{m-1})$ and $\cB=(B_{0},\ldots, B_{m-1})$. 
  Necessarily, we have that $\pm1\neq d \neq \pm 3$.
  Assume for a contradiction that $\cA$ and $\cB$ are equivalent, hence there is 
  a triple $(\varphi, g, z)$, where $\varphi$ is an automorphism of $G$, $g\in G$, and $z\in \Z_m$ such that
$\varphi(A_i)+g = B_{i+z}$, for every $i\in\Z_m$. Since $\pm 3$ is the only entry of the pattern appearing exactly once, it follows that $z=0$, hence
$\varphi(x) + g = \pm x$ for every $x\in\{1,3,d\}$. 
Assume that 
\begin{equation}\label{eq:lem:ineq:1}
\varphi(1) + g= (-1)^i
\;\;\;\text{and}\;\;\; 
\varphi(3) + g = (-1)^{i+j}3,
\end{equation} 
for some $i,j\in\{0,1\}$.
Since $\varphi(x)=x\varphi(1)$,  for every $x\in \Z_{m\ell^2+1}$, it follows that 
\[
2\varphi(1) = \varphi(3) + g - (\varphi(1) + g)
=(-1)^i (3(-1)^{j} - 1) = (-1)^{i+j}2^{j+1},
\]
hence $\varphi(1)=(-1)^{i+j}2^{j}$. Furthermore, by \eqref{eq:lem:ineq:1}, we have that
\begin{align*}
g = (-1)^i - \varphi(1)
\;\;\;\text{and}\;\;\;  
g =  (-1)^{i+j}3 - 3\varphi(1)
\end{align*}
hence,
\begin{align*}
2g&=3g-g = 3\big((-1)^i - \varphi(1)\big) - \big(3(-1)^{i+j} - 3\varphi(1) \big) \\
&= 3(-1)^i -3(-1)^{i+j}=3(-1)^i (1-(-1)^{j}) = 6(-1)^i j,
\end{align*}
therefore, $g=3(-1)^i j$. By recalling that $\cA\neq \pm\cB$, it follows that
$(\varphi,g)\neq (\pm id, 0)$, hence $j\neq 0$ (otherwise $g=0$, and $\varphi(1) = \pm 1$, that is, $\varphi=\pm id$). It then follows that $j=1$, hence 
\[
\text{$g=(-1)^i3$\;\;\; and\;\;\; $\varphi(1)=(-1)^{i+1}2$.}
\]
Finally, let $\varphi(d)+g=(-1)^kd$, that is, 
\[g=(-1)^kd - d\varphi(1) = d((-1)^k+(-1)^{i}2).\]
Since $g=(-1)^i3\neq \pm d$, then $i=k$, hence $(-1)^i3 = g=(-1)^i3d$, which implies $3(d-1)=0$, thus contradicting the assumption.
\end{proof}

Letting $\cA$ denote the CEDF built in Theorem \ref{main:0}, 
the following result provides a set $I$ of indices such that 
the sequence $\cA_i$ $i$-associated to $\cA$ 
consists of pairwise disjoint sets, for every $i\in I$; therefore,  by Lemma \ref{lem:associatedCEDF}, we obtain a family 
$\{\cA_i\mid i\in I\}$ of CEDFs with the same parameters as $\cA$. By checking that each $\cA_i$ satisfies the assumptions of Lemma \ref{lem:ineq:1}, it follows that 
$\{\cA\}\,\cup\,\{\cA_i\mid i\in I\}$ is a set of pairwise inequivalent CEDF.
This outlines the proof of the following result, which, due to its technical nature, is postponed to~\ref{Appendix:B}.

\begin{theorem} \label{th:inequivalent}
Let $m=24q+8y+z\geq 33$, with $y\in[0,2]$
   and $z\in[1,7]_2$. Also, let $\cA$ be the $(4m+1, m, 2, 1)$-CEDF built in Theorem \ref{main:0} and set 
   \[ 
   \textstyle
   I= [8q+4y+4, 16q+4y-4]_2\setminus\{12q+4y+2\lfloor\frac{z}{7}\rfloor\}.
   \]
   For each $i \in I$, let $\mathcal{A}_i$ be the sequence $i$-associated to $\mathcal{A}$.
   Then $\{\cA\}\,\cup\,\{\cA_i\mid i\in I\}$ is a set of $4q-3$  pairwise inequivalent 
   $(4m+1, m, 2, \lambda)$-CEDF.
\end{theorem}

\section{Conclusion}\label{sec:conclusions}
In this paper, we investigated the existence of cyclic $(v,m,\ell,1)$-CEDFs in the previously unresolved case where $m$ is odd and $\ell$ is even. Our main contribution is to provide explicit constructions in this setting, extending the known existence results.

More precisely, we proved that cyclic $(4m+1,m,2,1)$-CEDFs exist for every odd $m \geq 3$ 
(Theorems~\ref{main:1} and \ref{main:0}), 
and that cyclic $(3\ell^2+1,3,\ell,1)$-CEDFs exist for every even $\ell \geq 2$ 
(Theorem~\ref{main:2}). 
In both cases, the constructions are arithmetic and rely on a careful analysis of the interaction between the steps of the underlying arithmetic progressions. This step-based approach, already implicit in earlier works, is developed here as a systematic tool, leading in particular to the notions of \emph{pattern} and \emph{step-count}.

One outcome of this perspective is that it provides a simple invariant for distinguishing inequivalent CEDFs: arithmetic CEDFs with different step-counts are necessarily inequivalent. This allows us to exhibit multiple inequivalent families with the same parameters, and to systematically construct new ones by modifying suitable blocks. In particular, Theorem \ref{th:inequivalent} determines a subset $I\subset[0,m-1]$ such that ${\cal F}_I = \{\cA\}\,\cup\,\{\cA_i\mid i\in I\}$ is a family of pairwise inequivalent $(4m+1, m, 2, 1)$-CEDF, where $\cA$ is the CEDF built in Theorem \ref{main:0}, and $\cA_i$ is its $i$-associated sequence. This result shows that, even in the case $\ell=2$, the number of inequivalent cyclic CEDFs with fixed parameters can grow significantly with $m$. 
A larger set $I\subset [0,m-1]$ such that ${\cal F}_I$ is a family of pairwise inequivalent $(4m+1, m, 2, 1)$-CEDF is given in Table \ref{table}, for $m\geq 25$.

\begin{table}[h!]
\centering
\begin{tabular}{|c|l|}
\hline
$m$ & $I$ \\
\hline
$24q+1$ & $\{4q+1,\, 20q-1\} \cup ([8q+3,\,16q-3]_2 \setminus \{12q+1\})$ \\
\hline
$24q+3$ & $\{4q+1,\, 20q\pm1\} \cup ([8q+3,\,16q-3]_2 \setminus \{12q+1\})$ \\
\hline
$24q+5$ & $\{4q+1,\,12q+2\} \cup ([8q+3,\,16q-1]_2 \setminus \{12q+1\})$ \\
\hline
$24q+7$ & $\{4q+1,4q+3,\,20q+5,20q+3\} \cup ([8q+5,\,16q+1]_2 \setminus \{12q+3\})$ \\
\hline
$24q+9$ & $\{4q+3,\,20q+5\} \cup ([8q+5,\,16q+1]_2 \setminus \{12q+5\})$ \\
\hline
$24q+11$ & $\{4q+1,4q+3,\,20q+7\} \cup ([8q+5,\,16q+3]_2 \setminus \{12q+5\})$ \\
\hline
$24q+13$ & $\{4q+1,4q+3,12q+6,\,20q+7,20q+9\} \cup ([8q+7,\,16q+5]_2 \setminus \{12q+5\})$ \\
\hline
$24q+15$ & $\{4q+3\} \cup ([8q+7,\,16q+5]_2 \setminus \{12q+7\})$ \\
\hline
$24q+17$ & $\{4q+3,\,20q+13\} \cup ([8q+7,\,16q+7]_2 \setminus \{12q+9\})$ \\
\hline
$24q+19$ & $([8q+9,\,16q+9]_2 \setminus \{12q+9\})$ \\
\hline
$24q+21$ & $\{12q+10,\,20q+15\} \cup ([8q+9,\,16q+9]_2 \setminus \{12q+9\})$ \\
\hline
$24q+23$ & $\{20q+17\} \cup ([8q+9,\,16q+11]_2 \setminus \{12q+11\})$ \\
\hline
\end{tabular}
\caption{For every odd value of $m \pmod{24}$ ($m\geq 25$), this table provides
a subset $I\subset [0,m-1]$ such that ${\cal F}_I$ is a family of pairwise inequivalent $(4m+1, m, 2, 1)$-CEDF.
}\label{table}
\end{table}

As a final remark, every $(v,3,\ell,1)$-CEDF induces a
$(v,3,\ell,2)$-EDF. As a consequence, in Theorem \ref{EDF} we completely solve the existence problem for 
$(3\ell^2+1,3,\ell,2)$-EDFs over an abelian group.

Several natural questions remain open, in particular, the existence of
 a cyclic $(m\ell^2+1,m,\ell,1)$-CEDF for every odd $m>3$ and even $\ell>2$. 
While our results completely settle the cases $\ell=2$ and $m=3$, the general case remains elusive. 
For the same parameters, a stronger question is whether an arithmetic CEDF exists.

It could be also interesting to study
to what extent the method of modifying arithmetic CEDFs can be pushed further, for instance what  is the maximum number of pairwise inequivalent CEDFs that can be obtained from a single construction:
Section \ref{sec:inequivalent}
 provides a first answer in this direction.

The blow-up process described in \cite{PS2} (see also \cite{KPS} for an application to strong EDFs) can produce CEDFs that are not arithmetic, but whose sets are \emph{progressions of progressions}, that is, sets of size $\ell_1\cdots \ell_u>0$ of the following form
\begin{equation} \label{pro_of_pro}
[0,\ell_1-1]d_1+\cdots+[0,\ell_u-1]d_u+ x
\end{equation}
where the $d_i$s are nonzero integers and $x\in \Z$. 
It would therefore be interesting to investigate whether the notion of step count introduced here admits a suitable multidimensional generalization that applies to the broader class of CEDFs arising from this construction; in such a generalization, one might associate to a set as in \eqref{pro_of_pro} the corresponding sequence of steps $(d_1, \ldots, d_u)$. Such a generalization could potentially be exploited to obtain new infinite families of CEDFs extending those constructed in the present paper.

Recently, in \cite{HJM}, the authors consider digraph-defined EDF, where the differences are taken according to the arcs of a digraph $\Gamma$. In this setting a CEDF is a $\Gamma$-EDF, with $\Gamma$ a directed cycle. 
One might then investigate the relationship between the pattern of an arithmetic $\Gamma$-EDF and certain combinatorial properties of the associated graph $\Gamma$.
In particular, the chromatic number of $\Gamma$ provides a lower bound for the step-count. It would be interesting to determine whether a more precise relationship exists, and whether graph-theoretic invariants can be used to distinguish inequivalent EDFs or guide new constructions.

Overall, the results of this paper suggest that arithmetic structure and step patterns may play an important role in constructing CEDFs; this perspective offers a useful framework for further work on CEDFs and possibly also for digraph-defined EDFs.

\section{Acknowledgements}

The authors are grateful to the anonymous referee for the valuable comments and suggestions provided.

Burgess gratefully acknowledges support from NSERC Discovery Grant RGPIN-2025-04633.

Merola gratefully acknowledges support from project SERICS (PE00000014)
under the MUR National Recovery and Resilience Plan funded by the
European Union - NextGenerationEU.

Traetta gratefully acknowledges support from INdAM-GNSAGA.

\appendix

\section{Proof that the list $\mathcal{A}$ from Theorem~\ref{main:0} is a CEDF} \label{Appendix}

Here, we complete the proof of Theorem~\ref{main:0} for the case $m \geq 9$.  

Recall that $m=8b+2\epsilon+1$, where $\epsilon\in[0,3]$ and $b\geq 1$.  
For every $i \in \Z_m$, 
\[
x_i=
\begin{cases}
  2i & \text{if $i\in[0,m-3]_2$},\\
  2m-2i+6    & \text{if $i\in[3,4b-1]_4$},  \\
  2m-2i-10   & \text{if $i\in[1,4b-7]_4$},  \\  
  2m-8b-4\lceil \frac{\epsilon}{4}\rceil
  & \text{if $i=4b-3$},  \\
  4m-2i-10 
  & \text{if $\epsilon \neq 3$ and $i\in[4b+1, m-4]_2$} \\
	& \;\;\;\text{or if $\epsilon=3$ and $i \in [4b+5,m-4]_2$},\\
	m+7 & \text{if $\epsilon=3$ and $i=4b+1$}, \\
	3m-5 & \text{if $\epsilon=3$ and $i=4b+3$}, \\
   2m-7  & \text{if $i=m-2$},\\
   2m-1   & \text{if $i=m-1$}.  \\
\end{cases}
\]
Note that when $i=4b-3$ and $\epsilon \neq 0$, $x_i = 2m-2i-10$; if $\epsilon=0$, $x_{i} = 2m-2i-6$.
Defining $(d_0, d_1, \ldots, d_{m-1}) = (1,2m-2,1,\ldots,2m-2,1,3,2m-2)$, we have that $\mathcal{A} = (A_0, A_1, \ldots, A_{m-1})$, where 
for every $i \in \Z_m$, $A_i= x_i + \{0,d_i\}$.

We start by showing that the $A_i$s are pairwise disjoint. We partition the interval $[0,m-1]$ into the following four sets:
\[
I_0=[0,m-3]_2, \;I_1=[1,4b-1]_2, \; I_2=[4b+1, m-4]_2,\; I_3=\{m-2,m-1\},
\] 
and set $B_j= \bigcup_{i\in I_j} A_i$, for $j\in[0,3]$. 
One can easily check that
\[
\begin{aligned}
B_0 &= [0,2m-6]_4 \cup [1,2m-5]_4 = [0,16b+4\epsilon-4]_4 \cup [1,16b+4\epsilon-3]_4, \,\text{and}\\
B_3  &= \{2m-7, 2m-4, 2m-1, 4m-3\}.
\end{aligned}
\]
Let $X_j=\{x_i\mid i\in I_j\}$, for $j\in[1,2]$. 
By recalling that for $i\in [1,m-4]_2$, $A_i=\{x_i, x_i+d_i\}$ where $d_i= 2m-2$, it follows that $B_j= X_j\,\cup\, (X_j + 2m-2)$, 
for $j=1,2$. One can check that
\[
\begin{aligned}
  X_1 &=
  \textstyle{
  \{2m-8b-4\lceil\frac{\epsilon}{4}\rceil\}
    \,\cup\, \big([2m-8b+4,2m]_4\setminus\{2m-4\}\big),\;\text{and}
  }\\
  X_2 &=
  \begin{cases}
  [2m-2,4m-8b-12]_4 & \text{if $\epsilon \neq 3$}, \\
  [2m-2,4m-8b-20]_4 \cup \{m+7,3m-5\} & \text{if $\epsilon = 3$} 
  \end{cases}
   \\
  &=
    \begin{cases}
    [2m-2,4m-8b-12]_4 & \text{if $\epsilon \neq 3$}, \\
    [2m-2,4m-8b-20]_4 \cup \{m+7,4m-8b-12\} & \text{if $\epsilon = 3$} 
    \end{cases}
     \\
    &\subseteq [2m-2,4m-8b-12]_4 \cup
      \begin{cases}
      \varnothing & \text{if $\epsilon \neq 3$}, \\
      \{m+7\} & \text{if $\epsilon = 3$} 
      \end{cases}
       \\
\end{aligned}
\]
hence, $X_1+2m-2  \subset[4m-8b-6, 4m-2]_4$, and
\[
\begin{aligned}
X_2+2m-2 &  \subseteq [4m-4,6m-8b-14]_4 \cup 
\begin{cases}
\varnothing & \text{if $\epsilon \neq 3$} \\
\{3m+5\} & \text{if $\epsilon=3$}
\end{cases}
 \\
& = \{4m-4,4m\} \cup [3,2m-8b-15]_4 \cup
\begin{cases}
\varnothing & \text{if $\epsilon \neq 3$} \\
\{3m+5\} & \text{if $\epsilon = 3$}
\end{cases}
\\
\end{aligned}
\]
Since $b>0$, then $4m-8b-6> 2m$; therefore, $X_j$ and $X_j+2m-2$ are disjoint, hence $|B_j|= 2|X_j| = 2|I_j|$, for $j=1,2$. Also, the elements in $B_1$ are congruent to $2 \pmod{4}$; those in $B_2$ are congruent to $0$ or $3 \pmod{4}$, except for $m+7$ and $3m+5$ (when $m \equiv 7 \pmod{8}$) which are congruent to $2\pmod{4}$ 
 but do not belong to $B_1$.
Therefore, we have that $B_1$ and $B_2$ are disjoint. Noting that the elements in $B_0$ are congruent to $0$ or $1 \pmod{4}$, one can check that $B_0,B_1, B_2, B_3$ are pairwise disjoint.
Therefore,
\[
\left|\bigcup_{i=0}^{m-1} A_i\right| =  
\left|\bigcup_{j=0}^{3} B_j\right| = \sum_{j=0}^{3} |B_j| = \sum_{j=0}^{3} |I_j| =2m
\]
which implies that the $A_i$s are pairwise disjoint. 

It is left to show that $\bigcup_{i\in \Z_m}\Delta (A_{i}, A_{i-1}) = 
\Z_{4m+1}\setminus\{0\}$. 
Since
\[
\bigcup_{i=m-2}^{m} \Delta (A_{i}, A_{i-1}) 
= 
[1,6] \,\cup\, [2m+1, 2m+4] \,\cup\, [4m-1,4m],
\]
where $A_m=A_0$, it is enough to show that
\begin{equation}\label{eq:main:0}
\bigcup_{i=1}^{m-3}\Delta (A_{i}, A_{i-1}) 
= [7,2m]\,\cup\,[2m+5, 4m-2].
\end{equation}
Letting $T=\{0,1,2-2m, 3-2m\}$ and recalling that
$d_i=1$ or $2m-2$ according to whether $i\in [0,m-3]$ is even or odd,
one can check that
\[
\Delta(A_{i}, A_{i-1}) = (x_{i}-x_{i-1}) + 
\Delta(\{0, d_{i}\}, \{0, d_{i-1}\}) =
(x_{i}-x_{i-1}) +(-1)^i T.
\]

We first consider the case where
$m\not\equiv 7\pmod{8}$, that is, $\epsilon\in\{0,1,2\}$, and partition
$[1,m-3]$ into 
\[
J_0=[2,4b]_2,\,
J_1=[4b+2, m-3]_2,\;
J_2=[1,4b-1]_2,\;
J_3=[4b+1, m-4]_2.\;
\]
By taking into account that $T = (3-2m)-T$, we compute the following lists of differences:
\begin{align*}
\Delta_0:=&\;\textstyle{\bigcup_{i\in J_0}} \Delta (A_i, A_{i-1})
= \{x_i-x_{i-1} \mid i\in J_0\} + T = \\
=
&\; \big(\{2i-(2m-2(i-1)+6)\mid i\in[4,4b]_4\}    + T\big)  \,\cup\\
&\; \big(\{2i-(2m-2(i-1)-10)\mid i\in[2,4b-6]_4\} + T\big)   \,\cup\\
&\; \{2(4b-2)-(2m-8b \;-4
    \lceil
    \textstyle{\frac{\epsilon}{4}}
    \rceil)\} + T\big)  \\
=&\; \big(
     [-2m+8,-2m+16b-8]_8 
     \,\cup\, 
     \{-2m+16b+4\lceil
    \textstyle{\frac{\epsilon}{4}}
    \rceil-4\}
     \big) + T\\
=&\; \big(
     [2m+9, 2m+16b-7]_8 
     \,\cup\,      
     \{2m+16b+
     4
     \lceil
    \textstyle{\frac{\epsilon}{4}}
    \rceil-3\}
    \big)
     + T;\\
=&\; \big(
     [12, 16b-4]_8 
     \,\cup\,      
     \{16b+
     4
     \lceil
    \textstyle{\frac{\epsilon}{4}}
    \rceil\}
    \big)    
     - T;     \\ 
\Delta_1 :=&\;\textstyle{\bigcup_{i\in J_1}} \Delta (A_i, A_{i-1})
= \{x_i-x_{i-1} \mid i\in J_1\} + T = \\
= &\; \{2i-(4m-2(i-1)-10) \mid i \in J_1\} +T\\
= & \; \{4i-4m+8 \mid i \in J_1\} + T \\
=&\; \{4i+9\mid i\in[4b+2,m-3]_2\} + T\\
=&\; [16b+17,4m-3]_8 + T = [18-4\epsilon,2m]_8 - T;\\
\Delta_2:=&\;\textstyle{\bigcup_{i\in J_2}} \Delta (A_i, A_{i-1})
= \{x_i-x_{i-1} \mid i\in J_2\} - T = \\
=&\; \big(\{2m-2i+ 6 - 2(i-1) \mid i\in[3,4b-1]_4\} - T\big) \,\cup\\
&\; \big(\{2m-2i-10 - 2(i-1) \mid i\in[1,4b-7]_4\} - T\big) \,\cup\\
&\; \{(2m-8b -4
    \lceil
    \textstyle{\frac{\epsilon}{4}}
    \rceil) - 2(4b-4)
    \} + T\big)  \\
=&\;\big([2m-16b+12,2m-4]_8 \,\cup\,
    \{2m-16b+8-
    4
    \lceil
    \textstyle{\frac{\epsilon}{4}}
    \rceil
    \}
    \big) - T\\
=&\;
\big( 
[4\epsilon+14, 2m-4]_8 \,\cup\,
\{4\epsilon+10-
    4
    \lceil
    \textstyle{\frac{\epsilon}{4}}
    \rceil
\}
\big) - T\\
\Delta_3:=&\;\textstyle{\bigcup_{i\in J_3}} \Delta (A_i, A_{i-1})
= \{x_i-x_{i-1} \mid i\in J_3\} - T = \\
=&\; \{4m-2i-10-2(i-1) \mid i \in J_3\} -T\\
=&\; \{4m-8 -4i\mid i\in[4b+1,m-4]_2\} - T\\
=&\; [8, 4m -16b-12]_8 - T = [8, 16b+8\epsilon-8]-T.
\end{align*}
It is not difficult to check that
$\Delta_0\,\cup\, \Delta_3= [8,16b+4\epsilon]_4 -T =[8,2m-2]_4 - T$
and
$\Delta_1\,\cup\, \Delta_2 = [10,2m]_4-T$ 
; 
therefore,
\[
\begin{aligned}
\textstyle{
\bigcup_{i=1}^{m-3}}\Delta (A_{i}, A_{i-1}) 
&= 
\textstyle{\bigcup_{j=0}^{3}} \Delta_j = [8,2m]_2 - T\\
&=[8,2m]_2 + \{-1,0,2m-3, 2m-2\} \\
&= [7,2m]\,\cup\, [2m+5, 4m-2],
\end{aligned}
\]
thus proving \eqref{eq:main:0}. It is left to deal with the case where
$m\equiv 7\pmod{8}$, that is, $\epsilon=3$. We partition
$[1,m-3]$ into the following six subsets
\[
\begin{aligned}
J_0=[2,4b]_2,\;\;\;
J_1=\{4b+2, 4b+4\},\;\;\;
J_2=[4b+6, m-3]_2, \\
J_3=[1,4b-1]_2,\;\;\;
J_4=\{4b+1, 4b+3\},\;\;\;
J_5=[4b+5, m-4]_2,
\end{aligned}
\]
and recalling that $T = (3-2m)-T = (-16b-11)-T$ we compute the following lists of differences:
\begin{align*}
\Delta_0:=&\;\textstyle{\bigcup_{i\in J_0}} \Delta (A_i, A_{i-1})
= \{x_i-x_{i-1} \mid i\in J_0\} + T = \\
=
 &\; \big(\{2i-(2m-2(i-1)+6)\mid i\in[4,4b]_4\}    + T\big)  \,\cup\\
 &\; \big(\{2i-(2m-2(i-1)-10)\mid i\in[2,4b-2]_4\} + T\big)   \\
=&\; \big(
     [-2m+8,-2m+16b-8]_{16}\,\cup\, [-2m+16,-2m+16b]_{16}
     \big)+T\\
=&\; [-2m+8,-2m+16b]_8 + T = [12, 2m-10]_8 - T;\\
\Delta_1 :=&\;\textstyle{\bigcup_{i\in J_1}} \Delta (A_i, A_{i-1})
= \{x_{4b+2}-x_{4b+1}, x_{4b+4}-x_{4b+3}\} + T \\
=&\; \{8b+4-(8b+14), 8b+8-(24b+16)\} + T = \{-10, -16b -8\} + T\\
=&\; \{2m-6, 10\} - T;\\
\Delta_2 :=&\;\textstyle{\bigcup_{i\in J_2}} \Delta (A_i, A_{i-1})
= \{x_i-x_{i-1} \mid i\in J_2\} + T = \\
= &\; \{2i-(4m-2(i-1)-10) \mid i \in J_2\}+T \\
= &\; \{4i-4m+8 \mid i \in J_2\}+T \\
=&\; \{4i+9\mid i\in[4b+6,m-3]_2\} + T\\
=&\; [16b+33,4m-3]_8 + T = [22,2m]_8 - T;\\
\Delta_3:=&\;\textstyle{\bigcup_{i\in J_3}} \Delta (A_i, A_{i-1})
= \{x_i-x_{i-1} \mid i\in J_3\} - T = \\
=&\; \big(\{2m-2i-10 - 2(i-1) \mid i\in[1,4b-3]_4\} - T\big) \,\cup\\
 &\; \big(\{2m-2i+ 6 - 2(i-1) \mid i\in[3,4b-1]_4\} - T\big)  \\
=&\; \big([2m-16b+4, 2m-12]_{16}  \,\cup\, [2m-16b+12, 2m-4]_{16}\big) -T\\
=&\; [2m-16b+4, 2m-4]_8 - T \\
=&\; [18, 2m-4]_8-T\\
\Delta_4 :=&\;\textstyle{\bigcup_{i\in J_4}} \Delta (A_i, A_{i-1})
= \{x_{4b+1}-x_{4b}, x_{4b+3}-x_{4b+2}\} - T = \\
=&\; \{(8b+14)-8b, 24b+16-(8b+4)\} - T
= \{14, 16b +12\} - T\\
=&\; \{14,2m-2\} - T;\\
\Delta_5:=&\; \textstyle{\bigcup_{i\in J_5}} \Delta (A_i, A_{i-1})
= \{x_i-x_{i-1} \mid i\in J_5\} - T = \\
=&\; \{4m-2i-10-2(i-1) \mid i \in J_5\}-T\\
=&\; \{4m -4i - 8\mid i\in[4b+5,m-4]_2\} - T\\
=&\; [8, 4m -16b-28]_8 - T = [8, 2m-14]_8-T.
\end{align*}
It is not difficult to check that
$\Delta_0\,\cup\, \Delta_5= [8,2m-10]_4 -T$
and
$\Delta_2\,\cup\, \Delta_3 = [18,2m]_4-T$ 
; 
therefore,
\[
\begin{aligned}
&\; \textstyle{
\bigcup_{i=1}^{m-3}}\Delta (A_{i}, A_{i-1}) 
=
\textstyle{\bigcup_{j=0}^{5}} \Delta_j\\
&= 
\big([8,2m-10]_4 \,\cup\, [18,2m]_4 \,\cup\, \{10, 14,2m-6,2m-2\}\big)- T\\
&= \big([8,2m-2]_4 \,\cup\, [10,2m]_4\big)- T\\
&=[8,2m]_2 + \{-1,0,2m-3, 2m-2\} 
= [7,2m]\,\cup\, [2m+5, 4m-2],
\end{aligned}
\]
thus proving \eqref{eq:main:0}, and this completes the proof.

\section{Proof of Theorem \ref{th:inequivalent}} \label{Appendix:B}
Let $m=24q+8y+z\geq 33$, with $y\in[0,2]$
   and $z\in[1,7]_2$; hence, $4m+1 =  96q + 32y+ 4z + 1$.
   We start by recalling the definition of the $(4m+1, m, 2, 1)$-CEDF $\cA$ 
built in Theorem \ref{main:0} (with $b=3q+y$ and $\epsilon = \frac{z-1}{2}$).
More precisely, $\cA=(A_0, \ldots, A_{m-1})$, with $A_i=\{x_i, x_i+d_i\}$, where
$(d_0, \ldots, d_{m-1}) = (1,2m-2,1,\ldots,2m-2,1,3,2m-2)$ and
\[x_i=
\begin{cases}
  2i 
  & \text{if } i\in[0,m-3]_2,\\

  2m-2i+6
  & \text{if } i\in[3,\,12q+4y-1]_4,\\

  2m-2i-10
  & \text{if } i\in[1,\,12q+4y-7]_4,\\

  2m-2i-6
  -4\left\lceil \frac{z-1}{8} \right\rceil
  & \text{if } i=12q+4y-3,\\

  4m-10-2i
  & \text{if } i\in[12q+4y+1,\,m-4]_2, \\
	
	  2m-7  
  & \text{if } i=m-2,\\

  2m-1   
  & \text{if } i=m-1;\\

\end{cases}
\]
also, when $z=7$, we replace two of these values by setting 
\[x_i=
\begin{cases}
m+7  = 24q+8y+z + 7 &\text{if $i=12q+4y+1$,}\\
3m-5 = 72q + 24y + 3z - 5 &\text{if $i=12q+4y+3$.}
\end{cases}\]
Noticing that 
$I= [8q+4y+4, 16q+4y-4]_2\setminus\{12q+4y+2\lfloor\frac{z}{7}\rfloor\}$ is a set of even integers that are all less than $m-3$, for every $i\in I$ we have that
\[
x_i=2i, \;\;\;   d_i=1,
\;\;\;\text{and}\;\;\;
d_{i-1}= d_{i+1}=2m-2.
\]
It is enough to show that each $A_i+z_i = \{x_i+z_i, x_i+z_i+1\}$, with $i\in I$, is disjoint from $\bigcup_{i=0}^{m-1} A_i$, where
\begin{equation}
\begin{aligned}\label{eq:last}
z_i &=x_{i+1}+x_{i-1} - 2x_i + (d_{i-1} - d_i) \\
&=x_{i+1}+x_{i-1} -4i + 2m-3.
\end{aligned}
\end{equation}
Indeed, the sets of 
$\cA_i = (A_0, \ldots, A_{i-1}, A_i+z_i, A_{i+1}, \ldots, A_{m-1})$ (the sequence $i$-associated to $\cA$) will be pairwise disjoint, thus making $\cA_i$ a CEDF  having the same pattern as $\cA$, that is, $(\pm 1, \pm d,\ldots, \pm 1,  \pm3, \pm d)$, where $d=2m-2$. Notice that $3(d-1) = 6m-9 \neq 0$ in $\Z_{4m+1}$ and it is not difficult to check that $\cA_i\neq \cA_j$ for every distinct 
$i,j\in I$. Therefore, Lemmas~\ref{lem:associatedCEDF} and~\ref{lem:ineq:1} will guarantee that 
$\{\cA\}\,\cup\,\{\cA_i\mid i\in I\}$ is a set of pairwise inequivalent $(4m+1, m, 2, \lambda)$-CEDF.

Following the proof of Theorem \ref{main:0} (see~\ref{Appendix}), with $b=3q+y$ and $\epsilon = \frac{z-1}{2}$, 
we have that
$\bigcup_{i=0}^{m-1}A_i$ can be partitioned into four sets $B_0, \ldots, B_3$, where 
\begin{align*}
  B_0\,\cup\,B_3 &=
  [0,2m-6]_4 \,\cup\, [1,2m-5]_4
  \,\cup\,
  \{2m-7, 2m-4, 2m-1, 4m-3\}, \\
  &=[0,48q+16y+2z-6]_4 \,\cup\, [1,48q+16y+2z-5]_4  \,\cup\, \\
  &\;\;\;\;\, \{48q+16y+2z-7, 48q+16y+2z-4, 48q+16y+2z-1, 96q+32y+4z-3\} \\[2pt]
  B_1 &\subset  
  [2m-8b-4,2m]_4\,\cup\, [4m-8b-6, 4m-2]_4 \\
  &=[24q+8y+2z-4,\;48q+16y+2z]_4
    \;\cup \\
  &\;\;\;\;\,  
   [72q+24y+4z-6,\;96q+32y+4z-2]_4.\\[2pt]
  B_2 & \subseteq [2m-2,4m-8b-12]_4 \, \cup\, [3,2m-8b-15]_4 \,\cup\, \{4m-4,4m\} \,\cup\, \\
  &\;\;\;\;\,   \{m+7,3m+5\} \\
&= [48q+16y+2z-2,\;72q+24y+4z-12]_4
\,\cup\, [3,\;24q+8y+2z-15]_4 \,\cup\,
\\
&\;\;\;\;\, 
 \{96q+32y+4z-4,\;96q+32y+4z\} \,\cup\, B_2',      
\end{align*}
where $B_2' = \{24q+8y+z+7,\;72q+24y+3z+5\}$ if $z=7$, otherwise $B_2'=\varnothing$.

 Set $C_0=\bigcup_{i=0}^{m-1} A_i \,\cap\, [0, 4m]_4$ and $C_h=\bigcup_{i=0}^{m-1} A_i \,\cap\, [h, 4m-4+h]_4$, for $h\in[1,3]$; in other words, each $C_h$ is the subset of $[0,4m]$ containing exactly the elements of $\bigcup_{i=0}^{m-1} A_i$ congruent to $h$ modulo $4$, respectively. One can check that 
\begin{align*}
  C_0 \subseteq &\; 
 [0,\;48q+16y+2z-6]_4
 \,\cup\\
    &\; 
   [48q+16y+2z-2,\;72q+24y+4z-12]_4
 \,\cup \\ 
    &\; 
\{\; 96q+32y+4z-4,\; 96q+32y+4z\},\\
 C_2 \subseteq &\; 
B_1 \,\cup\, B_2' \,\cup\, \{48q+16y+2z-4\},
\;\;\;\text{and}  \\
 C_3  = &\;
 [3,\;24q+8y-2z+13]_4 \,\cup\, \{48q + 16y + 2z - 7\}.
 \end{align*}
We partition $I$ into the three subsets defined below:
\[
J_1=[8q+4y+4, 12q+4y - 2]_2,\;\;\;
J_2=[12q+4y+4, 16q+4y-4]_2,\;\;\;\text{and}\;\;\;
\]
\[\textstyle
J_3=\{12q+4y, 12q+4y+2\}\setminus\{12q+4y+2\lfloor\frac{z}{7}\rfloor\}.
\]
and consider the following three cases.\\

\textbf{Case 1}: Let $i\in J_1$ and note that
\[
x_{i-1}+x_{i+1}=
\begin{cases}
4m-4i & \text{if $i\in\{12q+4y-2,12q+4y-4\}$ and $z=1$},\\
4m-4i-4 & \text{otherwise}.
\end{cases}
\]
Therefore, by \eqref{eq:last} we have that
$z_i=2m - 8i - (4+\alpha)$, where $\alpha \in \{0,4\}$, and
\begin{align*}
 x_i+z_i &= 2m - 6i - (\alpha+4) = 6(m-i)-(\alpha+3) \in W \subset [3,4m-1]_4,\;\text{with}
\end{align*}
\[
W=[72q+ 24y+6z-\alpha+9, 96q + 24y+6z-\alpha-27]_{12}.
\]
Therefore, $x_i+z_i\equiv 3 \pmod{4}$ and since $W\,\cap\, C_3=\varnothing$, then
$x_i+z_i\not \in\bigcup_{i=0}^{m-1} A_i$. 
Also, $x_i+z_i + 1 \equiv 0\pmod{4}$ and since 
$(W+1)\,\cap\, C_0=\varnothing$, then $x_i+z_i + 1\not \in\bigcup_{i=0}^{m-1} A_i$; hence, 
$A_i+z_i$ is disjoint from $\bigcup_{i=0}^{m-1} A_i$.\\ 

\textbf{Case 2}: Let $i\in J_2$ and note that
\begin{align*}
x_{i-1}+x_{i+1} &= 8m-20-4i = 4m-21 - 4i. 
\end{align*}
Therefore,
$z_i= (4m-21 - 4i) -4i + 2m-3 = 6m -8i-24$ and
\begin{align*}
 x_i+z_i &= 6(m -i)-24 \in W \subset [2,4m-2]_4,\;\text{where}
\end{align*}
\[
W=[48q+24y+6z,\;72q+24y+6z-48]_{12}.
\]
Therefore, $x_i+z_i\equiv 2 \pmod{4}$ and since $W\,\cap\, C_2=\varnothing$, then
$x_i+z_i\not \in\bigcup_{i=0}^{m-1} A_i$. 
Also, $x_i+z_i + 1 \equiv 3\pmod{4}$ and since 
$(W+1)\,\cap\, C_3=\varnothing$, then $x_i+z_i + 1\not \in\bigcup_{i=0}^{m-1} A_i$; hence, 
$A_i+z_i$ is disjoint from $\bigcup_{i=0}^{m-1} A_i$. \\

\textbf{Case 3}: Let $i\in J_3$. If $z\neq 7$, then
$i=12q+4y+2$. Hence,  $x_{i-1}+x_{i+1} = 8m-20-4i$ and
and
$z_i = (8m-20-4i) -4i + 2m-3 = 6m -8i-24$. Therefore, 
\begin{align*}
 x_i+z_i &= 6m -6i-24 = 72q + 24y + 6z - 36 \equiv 2\pmod{4}.
\end{align*}
As before, since $x_i+z_i \not\in C_2$ and $x_i+z_i+1\not\in C_3$, then $A_i+z_i$ is disjoint from 
$\bigcup_{i=0}^{m-1} A_i$. 

If $z=7$, then $i=12q+4y$. Recall that
$x_{i-1} = 2m - 2(i-1)+6$ and $x_{i+1}=m+7$;
hence,  $x_{i-1}+x_{i+1} = 3m-2i+15$ and
and
$z_i = (3m-2i+15) -4i + 2m-3 = 5m-6i+12$. Therefore, 
\begin{align*}
 x_i+z_i &= 5m-4i+12 = 72q+24y+47 \equiv 3\pmod{4}.
\end{align*}
As before, since $x_i+z_i \not\in C_3$ and $x_i+z_i+1\not\in C_0$, then $A_i+z_i$ is disjoint from 
$\bigcup_{i=0}^{m-1} A_i$. Hence, $A_i+z_i$ is disjoint from $\bigcup_{i=0}^{m-1} A_i$.
\end{document}